\documentclass[10pt]{amsart} 

\usepackage{amsmath}
\usepackage{amssymb}

\usepackage{graphicx}
\usepackage{xcolor}
\usepackage{subcaption}
\usepackage{multirow}
\usepackage{enumitem}
\usepackage{mathrsfs}
\usepackage[rightcaption]{sidecap}
\usepackage{rotating}
\usepackage{tikz}

\newtheorem{theorem}{Theorem}[section]

\theoremstyle{definition}
\newtheorem{example}[theorem]{Example}
\theoremstyle{remark}

\newtheorem*{casei}{Case I}
\newtheorem*{caseii}{Case II}
\newtheorem*{caseiii}{Case III}

\newcommand{\bfu}{\mathbf u}
\newcommand{\bff}{\mathbf f}
\newcommand{\bfg}{\mathbf g}

\newcommand{\dx}{\mathrm d}
\newcommand{\FD}{\mathrm{FD}}
\newcommand{\JS}{\mathrm{JS}}
\newcommand{\Z}{\mathrm Z}
\newcommand{\ZA}{\mathrm{ZA}}
\newcommand{\cfl}{\text{CFL}}

\title[JS-type and Z-type weights for central-upwind WENO4 schemes]{JS-type and Z-type weights for fourth-order central-upwind weighted essentially non-oscillatory schemes}
\author{Jiaxi Gu}
\address{Department of Mathematics $\&$ POSTECH MINDS (Mathematical Institute for Data Science), Pohang University of Science and Technology, Pohang 37673, Korea}
\email{jiaxigu@postech.ac.kr}

\author{Xinjuan Chen}
\address{Department of Mathematics, College of Science, Jimei University, Xiamen, Fujian 361021, China}
\email{chenxinjuan@jmu.edu.cn}

\author{Kwanghyuk Park}
\address{Graduate School of Artificial Intelligence $\&$ POSTECH MINDS (Mathematical Institute for Data Science), Pohang University of Science and Technology, Pohang 37673, Korea}
\email{pkh0219@postech.ac.kr}

\author{Jae-Hun Jung}
\address{Department of Mathematics $\&$ POSTECH MINDS (Mathematical Institute for Data Science), Pohang University of Science and Technology, Pohang 37673, Korea}
\email{jung153@postech.ac.kr}

\makeatletter
\@namedef{subjclassname@2020}{\textup{}2020 Mathematics Subject Classification}
\makeatother

\begin{document}

\maketitle

\begin{abstract}
The central-upwind weighted essentially non-oscillatory (WENO) scheme introduces the downwind substencil to reconstruct the numerical flux, where the smoothness indicator for the downwind substencil is of critical importance in maintaining high order in smooth regions and preserving the essentially non-oscillatory behavior in shock capturing.
In this study, we design the smoothness indicator for the downwind substencil by simply summing up all local smoothness indicators and taking the average, which includes the regularity information of the whole stencil.
Accordingly the JS-type and Z-type nonlinear weights, based on simple local smoothness indicators, are developed for the fourth-order central-upwind WENO scheme.
The accuracy, robustness, and high-resolution properties of our proposed schemes are demonstrated in a variety of one- and
two-dimensional problems.
\end{abstract}

\section{Introduction} \label{sec:intro}
We consider the numerical solution of hyperbolic conservation laws
\begin{equation} \label{eq:1D_hyperbolic}
 \frac{\partial}{\partial t} u(x,t) + \frac{\partial}{\partial x} f(u(x,t)) = 0,
\end{equation} 
with the initial condition $u(x,0) = u_0(x)$.
The solution to \eqref{eq:1D_hyperbolic} may develop discontinuities within the finite time even for the smooth initial data.
It is challenging for the classical high-order finite difference discretizations to achieve high resolution in smooth regions, as well as maintain a sharp, essentially non-oscillatory transition around discontinuities.

Jiang and Shu \cite{Jiang} proposed the fifth-order WENO scheme (WENO5-JS) in the finite difference (FD) framework with the smoothness indicators, which measures the regularity over each substencil.
Then a convex combination of third-order numerical fluxes adapts either to the fifth-order numerical flux in smooth regions, or to a third-order numerical flux that eliminates the spurious oscillations in the presence of discontinuities.
However, the nonlinear weights $\omega_k^\JS$ in \cite{Jiang} fail in satisfying the sufficient condition for fifth-order accuracy at critical points, and there is also room for improvement in reducing numerical dissipation.
To increase the accuracy of the nonlinear weights, Henrick et al. \cite{Henrick} introduced mapping functions applied to $\omega_k^\JS$. 
The resulting nonlinear weights $\omega_k^M$ meet the requirement of fifth-order accuracy, while the mapped WENO scheme (WENO5-M) produces sharper approximations around discontinuities.
With a different approach, Borges et al. \cite{Borges} designed the global smoothness indicator and devised Z-type nonlinear weights $\omega_k^\Z$ with the WENO5-Z scheme, which further decreases the dissipation near discontinuities.
In \cite{HuWangAdams}, Hu et al. used the downwind substencil to reconstruct the numerical flux so that the scheme adapts to the central scheme with sixth-order accuracy in smooth regions.
When there is a discontinuity in the stencil, the contribution of the downwind substencil is disregarded, though not entirely dismissible, and it returns to the upwind scheme.
The adaption between the central and upwind schemes (WENO6-CU) is done by assigning the smoothness indicator of the whole stencil, which includes the regularity information, to the local smoothness indicator of the downwind substencil.

In this study, we adopt the notion of central-upwind WENO scheme \cite{HuWangAdams,Hu,HuAdams,Huang,LiuShenPengZhang,Wang,WangDonLiWang,Zhao} and develop the fourth-order WENO scheme by merging the downwind substencil into the 3-point upwind stencil of third-order WENO scheme, which is a 4-point stencil.
Unlike taking the smoothness indicator of this 4-point stencil for the local smoothness indicator of the downwind substencil, we simply average the three local smoothness indicator, which also contains the regularity information of the 4-point stencil.
Following the form of fifth-order JS-type and Z-type nonlinear weights, we design the corresponding fourth-order weights for the central-upwind WENO scheme. 

The outline of this paper is as follows.
The third-order WENO schemes, WENO3-JS and WENO3-Z, are described in Section \ref{sec:upwind_WENO}.
In Section \ref{sec:weights}, we introduce the novel JS-type and Z-type nonlinear weights for the fourth-order WENO schemes, and provide a detailed discussion about them.
We give the numerical evidence of the performed analysis in Section \ref{sec:ne} and present concluding remarks in Section \ref{sec:conclusion}.

\section{Upwind WENO scheme} \label{sec:upwind_WENO}
\begin{figure}[h!]
\centering
\includegraphics[width=0.8\textwidth]{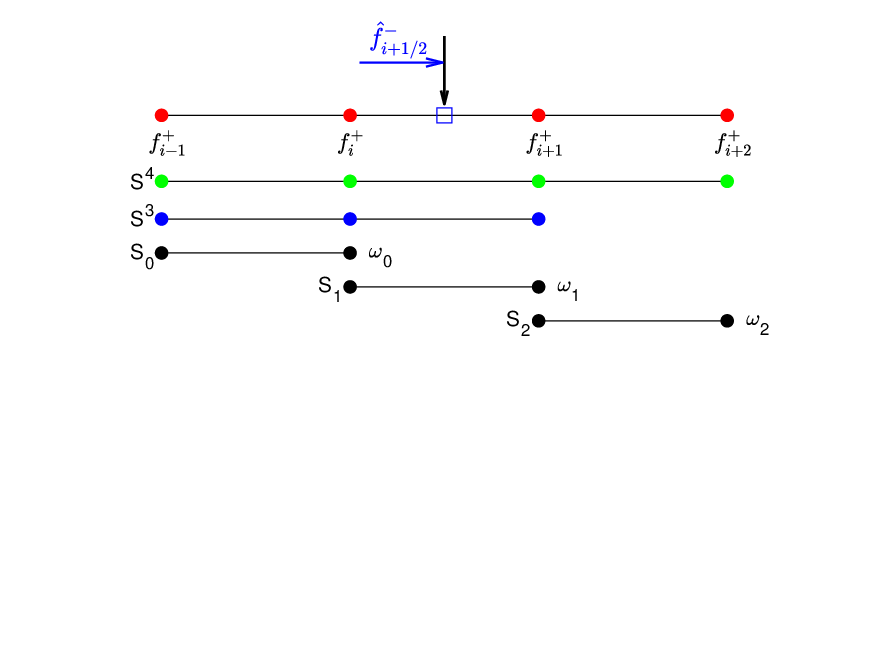}
\caption{The construction of numerical flux $\hat{f}^-_{i+1/2}$ depends on the stencils $S^3$ and $S^4$ for the respective third- and fourth-order accuracy, as well as three $2$-point substencils $S_0,\: S_1,\: S_2$ for second-order accuracy.}
\label{fig:stencil}
\end{figure}

We give a brief description of third-order WENO scheme to numerically solve \eqref{eq:1D_hyperbolic}.
Assuming that the time variable $t$ is continuous, we focus on the discussion of spatial discretization.
Given the computational domain $[a,\: b]$, we apply a uniform grid with $N$ points on the domain, 
$$
   x_i = a + \frac{1}{2} \Delta x + i \Delta x, \quad i = 1, \cdots, N, 
$$
where $\Delta x = \frac{b-a}{N}$ is the grid spacing. 
Each grid point $x_i$ is the cell center for the $i$th cell $I_i = \left[ x_{i-1/2},\: x_{i+1/2} \right]$ with the cell boundaries $x_{i \pm 1/2} = x_i \pm \Delta x/2$. 
If we fix the time $t$, the the method of lines discretization of \eqref{eq:1D_hyperbolic} yields
\begin{equation} \label{eq:1D_hyperbolic_discrete}
 \frac{du(x_i,t)}{dt} = - \left. \frac{\partial f \left( u(x,t) \right)}{\partial x} \right |_{x=x_i}.
\end{equation} 
Define the flux function $h(x)$ implicitly by
$$
   f \left( u(x) \right) = \frac{1}{\Delta x} \int^{x+\Delta x/2}_{x-\Delta x/2} h(\xi) \dx \xi,
$$
where we leave out $t$ in order to simplify the notation. 
Differentiating the flux function with respect to $x$ gives
$$
   \frac{\partial f}{\partial x} = \frac{h(x+\Delta x/2) - h(x-\Delta x/2)}{\Delta x},
$$
and \eqref{eq:1D_hyperbolic_discrete} becomes
\begin{equation} \label{eq:1D_hyperbolic_discrete_h}
 \frac{du(x_i,t)}{dt} = - \frac{h_{i+1/2} - h_{i-1/2}}{\Delta x},
\end{equation} 
with $h_{i \pm 1/2} = h(x_{i \pm 1/2})$. 
Since the fluxes $h_{i \pm 1/2}$ are defined in an implicit way, we aim to approximate explicitly.
Suppose that $\frac{\partial f}{\partial u} > 0$, i.e., the advection speed is positive.
A polynomial approximation $p^3(x)$ to $h(x)$ of degree at most two on the $3$-point stencil $S^3$, as shown in Figure \ref{fig:stencil}, is constructed.
Evaluating $p(x)$ at $x = x_{i+1/2}$ gives the finite difference numerical flux as
\begin{equation} \label{eq:numerical_flux_approximation_plus}
 \hat{f}^{\FD3}_{i+1/2} = - \frac{1}{6} f_{i-1} + \frac{5}{6} f_i + \frac{1}{3} f_{i+1},
\end{equation}
with $f_j = f(u_j)$.
The same procedure applied to the $2$-point substencils $S_0$ and $S_1$ in Figure \ref{fig:stencil} gives the numerical fluxes
\begin{equation} \label{eq:numerical_flux_approximation_plus_substencil}
\begin{aligned}
 \hat{f}^0_{i+1/2} &= - \frac{1}{2} f_{i-1} + \frac{3}{2} f_i, \\
 \hat{f}^1_{i+1/2} &=   \frac{1}{2} f_i     + \frac{1}{2} f_{i+1}.
\end{aligned}
\end{equation}
The convex linear combination of those low-order numerical fluxes with the linear weights $d_0 = \frac{1}{3}$ and $d_1 = \frac{2}{3}$ generates
$$ 
 \hat{f}^{\FD3}_{i+1/2} = d_0 \hat{f}^0_{i+1/2} + d_1 \hat{f}^1_{i+1/2}.
$$ 

The WENO numerical flux $\hat{f}_{i+1/2}$, which is an approximation of $h_{i+1/2}$, also takes a convex combination of two candidate numerical fluxes $\hat{f}^k_{i+1/2}$ in \eqref{eq:numerical_flux_approximation_plus_substencil}, 
\begin{equation} \label{eq:numerical_flux_approximation}
 \hat{f}_{i+1/2} = \omega_0 \hat{f}^0_{i+1/2} + \omega_1 \hat{f}^1_{i+1/2}.
\end{equation}
Note that the numerical flux $\hat{f}_{i+1/2}$ in \eqref{eq:numerical_flux_approximation} is $\hat{f}^-_{i+1/2}$ in Figure \ref{fig:stencil} as $\frac{\partial f}{\partial u} > 0$. 
For the stability and consistency, we require that
\begin{equation} \label{eq:weights}
 \omega_k \geqslant 0, \quad \omega_0 + \omega_1 = 1,
\end{equation}
where the sufficient condition for third-order accuracy in smooth regions is given by
\begin{equation} \label{eq:WENO3_condition}
 \omega_k - d_k = O(\Delta x^2).
\end{equation} 
In \cite{ShuSpringer}, the nonlinear weights $\omega^{\JS}_k$ in \eqref{eq:numerical_flux_approximation} are given by
\begin{equation} \label{eq:weights3_JS}
 \omega^{\JS}_k = \frac{\alpha_k}{\alpha_0 + \alpha_1},\quad \alpha_k = \frac{d_k}{(\beta_k + \epsilon)^2},\quad k=0,1,
\end{equation} 
where $\epsilon>0$ prevents the denominator from being zero, and the local smoothness indicators $\beta_k$ are
\begin{equation} \label{eq:SI}
 \beta_0 = \left( f_{i-1} - f_i \right)^2, \quad \beta_1 = \left( f_i - f_{i+1} \right)^2.
\end{equation}
Expanding \eqref{eq:SI} in Taylor series around $x=x_i$ gives
\begin{equation} \label{eq:SI_Taylor}
\begin{aligned}
 \beta_0 &= f'^2_i \Delta x^2 - f'_i f''_i \Delta x^3 + \left( \frac{1}{4} f''^2_i + \frac{1}{3} f'_i f'''_i \right) \Delta x^4 + O(\Delta x^5), \\
 \beta_1 &= f'^2_i \Delta x^2 + f'_i f''_i \Delta x^3 + \left( \frac{1}{4} f''^2_i + \frac{1}{3} f'_i f'''_i \right) \Delta x^4 + O(\Delta x^5).
\end{aligned}
\end{equation}
However, by Taylor expansion, $\omega^{\JS}_k = d_k + O(\Delta x)$, which does not satisfy the sufficient condition \eqref{eq:WENO3_condition}.
In \cite{Don}, Don and Borges introduced the global smoothness indicator $\tau_3$ as the absolute difference between $\beta_0$ and $\beta_1$,
$$
   \tau_3 = \left| \beta_0 - \beta_1 \right|,
$$
and defined the Z-type nonlinear weights $\omega^\Z_k$ as
\begin{equation} \label{eq:weights_Z}
 \omega^\Z_k = \frac{\alpha_k}{\alpha_0 + \alpha_1},\quad \alpha_k = d_k \left( 1 + \left( \frac{\tau_3}{\beta_k + \epsilon} \right)^2 \right),\quad k=0,1.
\end{equation}
Taylor series expansion shows that the weights $\omega^\Z_k$ satisfy $\omega^\Z_k = d_k + O(\Delta x^2)$, which agrees with the sufficient condition \eqref{eq:WENO3_condition}.

\section{JS-type and Z-type nonlinear weights for central-upwind WENO scheme} \label{sec:weights}
The central-upwind WENO scheme is constructed via adding the downwind substencil to the upwind stencil so that the scheme makes use of the central stencil for smooth regions while it reverts to the upwind stencil if there exists a discontinuity.
For the fourth-order central-upwind WENO scheme, we compute $p^4(x)$, the polynomial approximation to $h(x)$ over the 4-point stencil $S^4$ in Figure \ref{fig:stencil}, at $x = x_{i+1/2}$, which gives
\begin{equation} \label{eq:numerical_flux_approximation4_plus}
 \hat{f}^{\FD4}_{i+1/2} = - \frac{1}{12} f_{i-1} + \frac{7}{12} f_i + \frac{7}{12} f_{i+1} - \frac{1}{12} f_{i+2}.
\end{equation}
Applying the Taylor expansion, we obtain
\begin{align} \label{eq:numerical_flux_approximation4_plus_Taylor} 
 \hat{f}^{\FD4}_{i+1/2} = h_{i+1/2} - \frac{1}{30} h_i^{(4)} \Delta x^4 + O(\Delta x^5).
\end{align}
The numerical flux for the downwind substencil $S_2$ is  
\begin{equation} \label{eq:downwind_substencil}
 \hat{f}^2_{i+1/2} = \frac{3}{2} f_{i+1} - \frac{1}{2} f_{i+2},
\end{equation}
to \eqref{eq:numerical_flux_approximation_plus_substencil}.
Then we have the linear combination of the low-order numerical fluxes in \eqref{eq:numerical_flux_approximation_plus_substencil} and \eqref{eq:downwind_substencil} for the finite difference numerical flux \eqref{eq:numerical_flux_approximation4_plus}, such that
\begin{equation} \label{eq:numerical_flux_approximation4_linear_weights}
 \hat{f}^{\FD4}_{i+1/2} = \sum_{k=0}^2 d_k \hat{f}^k_{i+1/2},
\end{equation}
with the linear weights $d_0 = \frac{1}{6}, d_1 = \frac{2}{3}$ and $d_0 = \frac{1}{6}$.
We reconstruct the WENO numerical flux $\hat{f}_{i+1/2}$ of the form 
\begin{equation} \label{eq:numerical_flux_approximation4}
 \hat{f}_{i+1/2} = \sum_{k=0}^2 \omega_k \hat{f}^k_{i+1/2},
\end{equation}
where 
\begin{equation} \label{eq:weights4}
 \omega_k \geqslant 0,\quad \sum_{k=0}^2 \omega_k = 1.
\end{equation}

In the central-upwind WENO scheme, the smoothness indicator for the downwind substencil is critically important in maintaining high order in smooth regions and preserving the essentially non-oscillatory behavior in shock capturing.
Setting the downwind smoothness indicator as $\beta_d = \left( f_{i+1} - f_{i+2} \right)^2$, we define the smoothness indicator for the downwind substencil $\beta_2$ by
\begin{equation} \label{eq:SI_downwind}
 \beta_2 = \frac{1}{3} (\beta_0 + \beta_1 + \beta_d), 
\end{equation}
where the Taylor expansion around $x=x_i$ is
\begin{equation} \label{eq:SI_downwind_Taylor}
 \beta_2 = f'^2_i \Delta x^2 + f'_i f''_i \Delta x^3 + \left( \frac{11}{12} f''^2_i + f'_i f'''_i \right) \Delta x^4 + O(\Delta x^5).
\end{equation}

Following the form of the fifth-order JS-type nonlinear weights in \cite{Jiang}, we define the fourth-order JS-type weights as
\begin{equation} \label{eq:weights4_JS}
 \omega^{\JS}_k = \frac{\alpha_k}{\sum^2_{s=0} \alpha_s},\quad \alpha_k = \frac{d_k}{(\beta_k + \epsilon)^2},\quad k=0,1,2.
\end{equation}

Before designing the Z-type nonlinear weights, we first employ the methodology in \cite{ShuSpringer} to set the smoothness indicators $\beta_4$ for the stencil $S^4$ in Figure \ref{fig:stencil} as  
$$
   \beta_4 = \sum_{l=1}^3 \Delta x^{2l-1} \int_{x_{i-1/2}}^{x_{i+1/2}} \left( \frac{\dx^l}{\dx x^l} p^4_k(x) \right)^2 \dx x,
$$
which yields
\begin{equation} \label{eq:SI_high_order}
\begin{aligned}
 \beta_4 = & \frac{1}{9} \left( f_{i-1} - f_i - f_{i+1} + f_{i+2} \right)^2 + 
             \frac{44299}{103680} \left( f_{i-1} - 3 f_i + 3 f_{i+1} - f_{i+2} \right)^2 + \\
           & \frac{31}{57600} \left( f_{i-1} - 15 f_i + 15 f_{i+1} - f_{i+2} \right)^2 + 
              \frac{1}{2304} \left( 13 f_{i-1} + 29 f_i - 61 f_{i+1} + 19 f_{i+2} \right)^2 + \\
           & \frac{1}{2304} \left( 61 f_{i-1} - 151 f_i + 119 f_{i+1} - 29 f_{i+2} \right)^2 +         
             \frac{1}{32400} \left( 41 f_{i-1} -15 f_i + 15 f_{i+1} - 41 f_{i+2} \right)^2.
\end{aligned}
\end{equation}
Taylor series expansion of $\beta_4$ is
\begin{equation} \label{eq:SI_high_order_Taylor}
 \beta_4 = f'^2_i \Delta x^2 + \frac{13}{12} f''^2_i \Delta x^4 + O(\Delta x^5).
\end{equation}
We now define the global smoothness indicator $\tau_4$ on the stencil $S^4$ as
$$
   \tau_4 = \frac{1}{p} \left| \beta_4 - \frac{2 \beta_0 - 3 \beta_1 + 5 \beta_2}{4} \right|,
$$
and the new Z-type nonlinear weights $\omega^{\ZA}_k$  by
\begin{equation} \label{eq:weights4_ZA}
 \omega^{\ZA}_k = \frac{\alpha_k}{\sum^2_{s=0} \alpha_s},\quad \alpha_k = d_k \left( 1 + \left( \frac{\tau_4}{\beta_k + \epsilon} \right)^q \right),\quad k=0,1,2,
\end{equation}
with $p>0$ and $q \geqslant 1$. 

\subsection{Spatial fourth-order accuracy in smooth regions}
The finite difference numerical flux $\hat{f}^{\FD4}_{i-1/2}$ can be obtained by shifting one grid of $\hat{f}^{\FD4}_{i+1/2}$ \eqref{eq:numerical_flux_approximation4_plus} to the left
$$ 
   \hat{f}^{\FD4}_{i-1/2} = - \frac{1}{12} f_{i-2} + \frac{7}{12} f_{i-1} + \frac{7}{12} f_i - \frac{1}{12} f_{i+1},
$$ 
to which we apply the Taylor expansion,
\begin{equation} \label{eq:numerical_flux_approximation4_minus_Taylor}
 \hat{f}^{\FD4}_{i-1/2} = h_{i-1/2} - \frac{1}{60} h_i^{(5)} \Delta x^5 + O(\Delta x^6).
\end{equation}
By the similar analysis in \cite{Henrick}, we derive the sufficient condition for fourth-order accuracy in smooth regions as follows. 
Let
$$
   \hat{f}_{i \pm 1/2} = \sum_{k=0}^2 \omega^{\pm}_k \hat{f}^k_{i \pm 1/2}.
$$
Since
$$
   \hat{f}_{i+1/2} = \sum_{k=0}^2 d_k \hat{f}^k_{i+1/2} + \sum_{k=0}^2 ( \omega^+_k - d_k ) \hat{f}^k_{i+1/2} 
                   = \hat{f}^{\FD4}_{i+1/2} + \sum_{k=0}^2 ( \omega^+_k - d_k ) \hat{f}^k_{i+1/2},
$$
and
\begin{align*}
 \sum_{k=0}^2 ( \omega^+_k - d_k ) \hat{f}^k_{i+1/2} &= \sum_{k=0}^2 ( \omega^+_k - d_k ) \left[ h_{i+1/2} + A_k \Delta x^2 + O(\Delta x^3) \right] \\
                                                     &= h_{i+1/2} \sum_{k=0}^2 ( \omega^+_k - d_k ) + \Delta x^2 \sum_{k=0}^2 A_k ( \omega^+_k - d_k ) + \sum_{k=0}^2 ( \omega^+_k - d_k ) O(\Delta x^3) \\
                                                     &= \Delta x^2 \sum_{k=0}^2 A_k ( \omega^+_k - d_k ) + \sum_{k=0}^2 ( \omega^+_k - d_k ) O(\Delta x^3), 
\end{align*}
then
$$
   \hat{f}_{i+1/2} = \hat{f}^{\FD4}_{i+1/2} + \Delta x^2 \sum_{k=0}^2 A_k ( \omega^+_k - d_k ) + \sum_{k=0}^2 ( \omega^+_k - d_k ) O(\Delta x^3). 
$$
Similarly,
$$
   \hat{f}_{i-1/2} = \hat{f}^{\FD4}_{i-1/2} + \Delta x^2 \sum_{k=0}^2 A_k ( \omega^-_k - d_k ) + \sum_{k=0}^2 ( \omega^-_k - d_k ) O(\Delta x^3). 
$$
By \eqref{eq:numerical_flux_approximation4_plus_Taylor} and \eqref{eq:numerical_flux_approximation4_minus_Taylor}, we have
\begin{align*}
 \frac{\hat{f}_{i+1/2} - \hat{f}_{i-1/2}}{\Delta x} = {} & \frac{h_{i+1/2} - h_{i-1/2}}{\Delta x} + O(\Delta x^4) + \Delta x \sum_{k=0}^2 A_k ( \omega^+_k - \omega^-_k ) \\
                                                      {} & + \sum_{k=0}^2 ( \omega^+_k - d_k ) O(\Delta x^2) - \sum_{k=0}^2 ( \omega^-_k - d_k ) O(\Delta x^2).
\end{align*} 
Thus the sufficient condition for fourth-order accuracy is given by
\begin{equation} \label{eq:WENO4_condition}
 \omega_k - d_k = O(\Delta x^3),
\end{equation} 
where the superscripts are dropped, meaning that the nonlinear weights $\omega_k$ for each stencil $S^4$ should satisfy the condition \eqref{eq:WENO4_condition} in order to attain fourth-order accuracy in space.

For the nonlinear weights $\omega^{\JS}_k$ in \eqref{eq:weights4_JS}, it can be shown via Taylor expansion that $\omega^{\JS}_k = d_k + O(\Delta x)$, and thus the weights $\omega^{\JS}_k$ do not satisfy the sufficient condition \eqref{eq:WENO4_condition}.

Expanding the global smoothness indicator $\tau_4$ in a Taylor series at $x=x_i$ yields
$$
   \tau_4 = - \frac{7}{6} f'_i f'''_i \Delta x^4 - \left( \frac{5}{4} f''_i f'''_i + \frac{7}{12} f'_i f^{(4)}_i \right) \Delta x^5 + O(\Delta x^6).
$$
Combing this with \eqref{eq:SI_Taylor} and \eqref{eq:SI_downwind_Taylor}, we obtain $\omega^\ZA_k = d_k + O(\Delta x^{2q})$ if $f'_i \ne 0$, whereas $\omega^\ZA_k = d_k + O(\Delta x^q)$ if $f'_i = 0$.
In order to ensure that the weights $\omega^\ZA_k$ satisfy the condition \eqref{eq:WENO4_condition} at least at non-critical points, $q \geqslant \frac{3}{2}$.
As explained in \cite{ChenGuJung}, a more direct way to achieve fourth-order accuracy in smooth regions is to increase $p$ as much as possible.

\subsection{Weights behavior around discontinuity}
According to the WENO idea, when there is a discontinuity in one or more of the substencils, we expect the corresponding weight(s) to be essentially zero.
In \cite{Borges,Gu}, it was observed that if the WENO scheme assigns orders of magnitude larger weight(s) to discontinuous stencil(s), the convex combination of the low-order numerical fluxes is closer to the traditional finite difference numerical flux, leading to less dissipation near the discontinuity.
We also analyzed the effect of the WENO weighting procedure on the numerical dissipation in \cite{ChenGuJung}. 

To gain some insight into the behavior of the nonlinear weights $\omega_k$ for the numerical dissipation when a discontinuity exists, assume that the the solution is smooth over the stencil $S^4$ except a discontinuity point in the substencil $S_k$. 
Then there are three cases in our analysis.
By setting $\epsilon = 0$, \eqref{eq:weights4_JS} implies that
\begin{align}
 \omega^{\JS}_0 &= \frac{1}{d_0 + d_1 \left( \frac{\beta_0}{\beta_1} \right)^2 + d_2 \left( \frac{\beta_0}{\beta_2} \right)^2} d_0, \label{eq:weights4_JS_0} \\
 \omega^{\JS}_1 &= \frac{1}{d_0 \left( \frac{\beta_1}{\beta_0} \right)^2 + d_1 + d_2 \left( \frac{\beta_1}{\beta_2} \right)^2} d_1, \label{eq:weights4_JS_1} \\
 \omega^{\JS}_2 &= \frac{1}{d_0 \left( \frac{\beta_2}{\beta_0} \right)^2 + d_1 \left( \frac{\beta_2}{\beta_1} \right)^2 + d_2} d_2. \label{eq:weights4_JS_2}
\end{align}
Similarly, with $q=2$ for $\omega^\ZA_k$ \eqref{eq:weights4_ZA}, satisfying the sufficient condition for fourth-order accuracy \eqref{eq:WENO4_condition}, we have
\begin{align}
 \omega^\ZA_0 &= \frac{1}{d_0 + d_1 \frac{1 + \left( \frac{\tau_4}{\beta_1} \right)^2}{1 + \left( \frac{\tau_4}{\beta_0} \right)^2} + d_2 \frac{1 + \left( \frac{\tau_4}{\beta_2} \right)^2}{1 + \left( \frac{\tau_4}{\beta_0} \right)^2}} d_0, \label{eq:weights4_ZA_0} \\
 \omega^\ZA_1 &= \frac{1}{d_0 \frac{1 + \left( \frac{\tau_4}{\beta_0} \right)^2}{1 + \left( \frac{\tau_4}{\beta_1} \right)^2} + d_1 + d_2 \frac{1 + \left( \frac{\tau_4}{\beta_2} \right)^2}{1 + \left( \frac{\tau_4}{\beta_1} \right)^2}} d_1, \label{eq:weights4_ZA_1} \\
 \omega^\ZA_2 &= \frac{1}{d_0 \frac{1 + \left( \frac{\tau_4}{\beta_0} \right)^2}{1 + \left( \frac{\tau_4}{\beta_2} \right)^2} + d_1 \frac{1 + \left( \frac{\tau_4}{\beta_1} \right)^2}{1 + \left( \frac{\tau_4}{\beta_2} \right)^2} + d_2} d_2, \label{eq:weights4_ZA_2}
\end{align}
plus the additional assumption $\tau_4 = O(1)$.

\begin{casei}
We assume the discontinuity lies in the substencil $S_2$.
Then the smoothness indicators satisfy $\beta_0,\: \beta_1 = O(\Delta x^2)$ and $\beta_d,\: \beta_2 = O(1)$, from which we have $\tau_4 \approx \frac{167}{240p} \beta_d$ and 
\begin{equation} \label{eq:betak2_ineq}
 \beta_k \ll \beta_2,\quad 1 + \left(\frac{\tau_4}{\beta_k}\right)^2 \gg 1 + \left( \frac{\tau_4}{\beta_2}\right)^2,\quad k = 0,1.
\end{equation}
It follows that
\begin{equation} \label{eq:betak2}
 \frac{\beta_k}{\beta_2} = o(1),\quad \frac{1 + \left(\frac{\tau_4}{\beta_2}\right)^2}{1 + \left(\frac{\tau_4}{\beta_k}\right)^2} = o(1),\quad k = 0,1.
\end{equation}
Since $\displaystyle{\lim_{\Delta x \to 0} \beta_0 / \beta_1 = 1}$, we see that 
$$
   \lim_{\Delta x \to 0} \frac{1 + \left(\frac{\tau_4}{\beta_1}\right)^2}{1 + \left(\frac{\tau_4}{\beta_0}\right)^2} - 1
 = \lim_{\Delta x \to 0} \frac{\left[ 1 + \left(\frac{\tau_4}{\beta_1}\right)^2 \right] \left(\frac{\beta_0}{\tau_4}\right)^2}{\left[1 + \left(\frac{\tau_4}{\beta_0}\right)^2 \right] \left(\frac{\beta_0}{\tau_4}\right)^2} - 1
 = \lim_{\Delta x \to 0} \frac{\left(\frac{\beta_0}{\beta_1}\right)^2 - 1}{1 + \left(\frac{\beta_0}{\tau_4}\right)^2} = 0,
$$
which gives
\begin{equation} \label{eq:beta01}
 \frac{\beta_0}{\beta_1} = 1 + o(1),\quad \frac{1 + \left(\frac{\tau_4}{\beta_1}\right)^2}{1 + \left(\frac{\tau_4}{\beta_0}\right)^2} = 1 + o(1).
\end{equation}
With \eqref{eq:betak2} and \eqref{eq:beta01}, we find
$$
   d_0 + d_1 \left( \frac{\beta_0}{\beta_1} \right)^2 + d_2 \left( \frac{\beta_0}{\beta_2} \right)^2 = d_0 + d_1 \left( 1 + o(1) \right)^2 + d_2 \, o(1) = d_0 + d_1 + o(1) < d_0 + d_1 + d_2 = 1,
$$
and
$$
   d_0 + d_1 \frac{1 + \left( \frac{\tau_4}{\beta_1} \right)^2}{1 + \left( \frac{\tau_4}{\beta_0} \right)^2} + d_2 \frac{1 + \left( \frac{\tau_4}{\beta_2} \right)^2}{1 + \left( \frac{\tau_4}{\beta_0} \right)^2} = d_0 + d_1 \left( 1 + o(1) \right) + d_2 \, o(1) = d_0 + d_1 + o(1) < d_0 + d_1 + d_2 = 1.
$$
Consequently, by \eqref{eq:weights4_JS_0} and \eqref{eq:weights4_ZA_0},
\begin{equation} \label{eq:weights4_0_approx}
 \omega^\JS_0,~\omega^\ZA_0 > d_0,\quad \omega^\JS_0,~\omega^\ZA_0 = \frac{1}{5} + o(1).
\end{equation}
Using the same approach, we obtain 
\begin{equation} \label{eq:weights4_1_approx}
 \omega^\JS_1,~\omega^\ZA_1 > d_1,\quad \omega^\JS_1,~\omega^\ZA_1 = \frac{4}{5} + o(1).
\end{equation}
Note that for $k=0,1$,
$$
   \lim_{\Delta x \to 0} \frac{\beta_k^2 + \tau_4^2}{\beta_2^2 + \tau_4^2} = \frac{27889}{34289}.
$$
We thus expect that if $\Delta x$ is small enough,
\begin{equation} \label{eq:betak2_tau4_approx}
 \frac{\beta_k^2 + \tau_4^2}{\beta_2^2 + \tau_4^2} \approx \frac{27889}{34289}.  
\end{equation}
From $\beta_2 \gg \beta_k$ and \eqref{eq:betak2_tau4_approx}, we have
$$
   \left(\frac{\beta_2}{\beta_k}\right)^2 \left( 1 - \frac{\beta_k^2 + \tau_4^2}{\beta_2^2 + \tau_4^2} \right) \gg 1 > 0,
$$
and hence
$$
   \left( \frac{\beta_2}{\beta_k} \right)^2 \gg \left( \frac{\beta_2}{\beta_k} \right)^2 \frac{\beta_k^2 + \tau_4^2}{\beta_2^2 + \tau_4^2} = \frac{1 + \left(\frac{\tau_4}{\beta_k}\right)^2}{1 + \left(\frac{\tau_4}{\beta_2}\right)^2}.
$$
The above inequality and \eqref{eq:betak2_ineq} give
$$
   d_0 \left( \frac{\beta_2}{\beta_0} \right)^2 + d_1 \left( \frac{\beta_2}{\beta_1} \right)^2 + d_2 \gg 
   d_0 \frac{1 + \left( \frac{\tau_4}{\beta_0} \right)^2}{1 + \left( \frac{\tau_4}{\beta_2} \right)^2} + d_1 \frac{1 + \left( \frac{\tau_4}{\beta_1} \right)^2}{1 + \left( \frac{\tau_4}{\beta_2} \right)^2} + d_2 \gg 
   d_0 + d_1 + d_2 = 1,
$$
so $\omega^\JS_2 \ll \omega^\ZA_2 \ll d_2$ by \eqref{eq:weights4_JS_2} and \eqref{eq:weights4_ZA_2}.
Using the requirement \eqref{eq:weights4} and the estimates \eqref{eq:weights4_0_approx}, \eqref{eq:weights4_1_approx}, we get $\omega^\JS_2, \omega^\ZA_2 \gtrapprox 0$.

In summary,
\begin{gather*}
 \omega^\JS_0,~\omega^\ZA_0 \approx \frac{1}{5},\quad 
 \omega^\JS_1,~\omega^\ZA_1 \approx \frac{4}{5},\quad 
 0 \lessapprox \omega^\JS_2 \ll \omega^\ZA_2 \ll d_2.
\end{gather*}
\end{casei}

\begin{caseii}
Suppose that we have the discontinuity in the substencil $S_1$.
Since $\beta_0, \beta_d= O(\Delta x^2)$ and $\beta_1, \beta_2 = O(1)$, then
$$
   \beta_0 \ll \beta_k,\quad 1 + \left(\frac{\tau_4}{\beta_0}\right)^2 \gg 1 + \left(\frac{\tau_4}{\beta_k}\right)^2,\quad k = 1,2,
$$
which implies that 
\begin{equation} \label{eq:beta0k}
 \frac{\beta_0}{\beta_k} = o(1),\quad \frac{1 + \left(\frac{\tau_4}{\beta_k}\right)^2}{1 + \left(\frac{\tau_4}{\beta_0}\right)^2} = o(1),\quad k = 1,2.
\end{equation}
Applying \eqref{eq:beta0k} to $\omega^\JS_0$ in \eqref{eq:weights4_JS_0} and $\omega^\ZA_0$ in \eqref{eq:weights4_ZA_0}, we find
\begin{align*}
 \omega^\JS_0 &= \frac{1}{d_0 + d_1 \, o(1) + d_2 \, o(1)} d_0 = 1 + o(1), \\
 \omega^\ZA_0 &= \frac{1}{d_0 + d_1 \, o(1) + d_2 \, o(1)} d_0 = 1 + o(1).
\end{align*}
Combining this with the requirement \eqref{eq:weights4} gives $\omega^\JS_0, \omega^\ZA_0 \lessapprox 1$, and thus $\omega^\JS_k, \omega^\ZA_k \gtrapprox 0,\: k=1,2$.
Since, for $k=1,2$,
$$
   \left(\frac{\beta_0}{\beta_k}\right)^2 \left( 1 - \frac{\beta_k^2 + \tau_4^2}{\beta_0^2 + \tau_4^2}\right) 
 = \left(\frac{\beta_0}{\beta_k}\right)^2 \frac{\beta_0^2 - \beta_k^2}{\beta_0^2 + \tau_4^2} < 0,
$$
we obtain
$$
   \left(\frac{\beta_0}{\beta_k}\right)^2 < \left(\frac{\beta_0}{\beta_k}\right)^2 \frac{\beta_k^2 + \tau_4^2}{\beta_0^2 + \tau_4^2} = \frac{1 + \left(\frac{\tau_4}{\beta_k}\right)^2}{1 + \left(\frac{\tau_4}{\beta_0}\right)^2},
$$
and
$$
   d_0 + d_1 \left(\frac{\beta_0}{\beta_1}\right)^2 + d_2 \left(\frac{\beta_0}{\beta_2}\right)^2 
 < d_0 + d_1 \frac{1 + \left(\frac{\tau_4}{\beta_1}\right)^2}{1 + \left(\frac{\tau_4}{\beta_0}\right)^2} + d_2 \frac{1 + \left(\frac{\tau_5}{\beta_2}\right)^2}{1 + \left(\frac{\tau_5}{\beta_0}\right)^2} 
 < d_0 + d_1 + d_2 = 1,
$$
and hence $\omega^\JS_0 > \omega^\ZA_0 > d_0$ according to \eqref{eq:weights4_JS_0} and \eqref{eq:weights4_ZA_0}.
Now, with $\beta_1 \gg \beta_0$, we have
$$
   \left(\frac{\beta_1}{\beta_0}\right)^2 \left( 1 - \frac{\beta_0^2 + \tau_4^2}{\beta_1^2 + \tau_4^2} \right) \gg 
   \frac{d_2}{d_0} \left(\frac{\beta_1}{\beta_2}\right)^2 \left( \frac{\beta_2^2 + \tau_4^2}{\beta_1^2 + \tau_4^2} - 1 \right),
$$
and hence
\begin{equation} \label{eq:omega1_JS_ZA}
 d_0 \left( \frac{\beta_1}{\beta_0} \right)^2 + d_1 + d_2 \left( \frac{\beta_1}{\beta_2} \right)^2 \gg 
 d_0 \frac{1 + \left(\frac{\tau_4}{\beta_0}\right)^2}{1 + \left(\frac{\tau_4}{\beta_1}\right)^2} + d_1 + d_2 \frac{1 + \left(\frac{\tau_4}{\beta_2}\right)^2}{1 + \left(\frac{\tau_4}{\beta_1}\right)^2} \gg
 d_0 + d_1 + d_2 = 1.
\end{equation}
So $\omega^\JS_1 \ll \omega^\ZA_1 \ll d_1$.
A similar argument leads to the inequality $\omega^\JS_2 \ll \omega^\ZA_2 \ll d_2$.

In summary, we have
$$
 d_0 < \omega^\ZA_0 < \omega^\JS_0 \lessapprox 1,\quad
 0 \lessapprox \omega^\JS_k \ll \omega^\ZA_k \ll d_k,\ k = 1,2.
$$
\end{caseii}

\begin{caseiii}
Assume that the substencil $S_0$ contains the discontinuity.
Using the analysis that is analogous to what is done in Case II, we obtain the following:
$$
 d_1 < \omega^\ZA_1 < \omega^\JS_1 \lessapprox 1,\quad
 0 \lessapprox \omega^\JS_k \ll \omega^\ZA_k \ll d_k,\ k = 0,2.
$$
\end{caseiii}

As our above analysis indicates, the Z-type weights $\omega^\ZA_k$ are closer to the linear weights $d_k$ than the JS-type ones $\omega^\JS_k$.
Then the former convex combination is closer to the classical finite difference scheme than the latter, and thus we expect the sharper approximation near discontinuities with $\omega^\ZA_k$, which is verified in the next section.

\section{Numerical examples} \label{sec:ne}
We present one- and two-dimensional numerical results to demonstrate the potential of the proposed fourth-order central-upwind WENO schemes, denoted as WENO4-JS and WENO4-ZA.
We take $q=2$ for the nonlinear weights $\omega^\ZA_k$ in \eqref{eq:weights4_ZA} as the sufficient condition \eqref{eq:WENO4_condition} is satisfied for at least non-critical points. 
We use the one-dimensional linear advection equation to check the order of accuracy of the WENO schemes in terms of $L_1$ and $L_\infty$ error norms, where different values of $p$ and $\epsilon$ are set for WENO4-ZA.
The other examples show the numerical results from WENO4-JS and WENO4-ZA, in comparison with WENO3-Z and WENO5-JS.
We choose $\epsilon = 10^{-40}$ for WENO3-Z and WENO4-ZA whereas $\epsilon = 10^{-6}$ remains for WENO4-JS and WENO5-JS.
We integrate in time using the explicit third-order total variation diminishing Runge-Kutta method \cite{ShuOsherI} and set $\cfl = 0.4$.

\subsection{1D linear advection equation} 
\begin{example} \label{ex:advection_1d}
Our first example checks the order of accuracy for the one-dimensional linear advection equation
$$
   u_t + u_x = 0,\quad -1 \leqslant x \leqslant 1,
$$
subject to the initial condition $u(x,0) = \sin(\pi x)$ and the periodic boundary condition. 
The exact solution is given by $u(x,t)= \sin \left( \pi (x-t) \right)$.
The numerical solution is computed up to the time $T=2$ with the time step $\Delta t = \cfl \cdot \Delta x^{4/3}$. 
The $L_1$ and $L_\infty$ errors with the numerical order, for the WENO4-JS, WENO4-ZA1 ($p=$ 1e+2, $\epsilon=$ 1e-40), WENO4-ZA2 ($p=$ 1e+5, $\epsilon=$ 1e-16), and the classical fourth-order finite difference (FD4) schemes, are listed in Tables \ref{tab:advection_1d_L1} and \ref{tab:advection_1d_Linf}, respectively.
We see that WENO4-ZA and FD4 achieve the expected order of accuracy.
It also shows that the parameter $\epsilon$ plays a role affecting the numerical order of WENO schemes, as in \cite{Don}.
\end{example}

\begin{table}[h!]
\renewcommand{\arraystretch}{1.1}
\scriptsize
\centering
\caption{$L_1$ error and numerical order for Example \ref{ex:advection_1d}.}      
\begin{tabular}{clcrlcrlcrlc} 
\hline  
N &\multicolumn{2}{l}{WENO4-JS} & &\multicolumn{2}{l}{WENO4-ZA1} & &\multicolumn{2}{l}{WENO4-ZA2} & &\multicolumn{2}{l}{FD4} \\ 
   \cline{2-3}                     \cline{5-6}                      \cline{8-9}                      \cline{11-12}  
  & Error & Order               & & Error & Order                & & Error & Order                & & Error & Order \\
\hline
10  & 2.42e-1 & --     & & 5.23e-2 & --     & & 1.92e-2 & --     & & 1.93e-2 & -- \\  
20  & 9.66e-2 & 1.3259 & & 6.91e-3 & 2.9196 & & 1.30e-3 & 3.8892 & & 1.29e-3 & 3.9099 \\  
40  & 3.70e-2 & 1.3861 & & 3.28e-3 & 1.0736 & & 3.95e-4 & 1.7127 & & 8.10e-5 & 3.9906 \\
80  & 1.03e-2 & 1.8396 & & 8.15e-4 & 2.0100 & & 2.13e-5 & 4.2113 & & 5.07e-6 & 3.9975 \\ 
160 & 2.49e-3 & 2.0539 & & 2.00e-4 & 2.0298 & & 1.23e-6 & 4.1207 & & 3.17e-7 & 3.9990 \\ 
\hline
\end{tabular}
\label{tab:advection_1d_L1}
\end{table}

\begin{table}[h!]
\renewcommand{\arraystretch}{1.1}
\scriptsize
\centering
\caption{$L_\infty$ error and numerical order for Example \ref{ex:advection_1d}.}      
\begin{tabular}{clcrlcrlcrlc} 
\hline  
N &\multicolumn{2}{l}{WENO4-JS} & &\multicolumn{2}{l}{WENO4-ZA1} & &\multicolumn{2}{l}{WENO4-ZA2} & &\multicolumn{2}{l}{FD4} \\ 
   \cline{2-3}                     \cline{5-6}                      \cline{8-9}                      \cline{11-12}  
  & Error & Order               & & Error & Order                & & Error & Order                & & Error & Order \\
\hline
10  & 4.75e-1 & --     & & 1.12e-1 & --     & & 2.99e-2 & --     & & 2.99e-2 & -- \\  
20  & 1.91e-1 & 1.3121 & & 1.77e-2 & 2.6559 & & 2.86e-3 & 3.3864 & & 2.00e-3 & 3.9038 \\ 
40  & 7.98e-2 & 1.2627 & & 1.03e-2 & 0.7780 & & 1.22e-3 & 1.2253 & & 1.27e-4 & 3.9759 \\
80  & 3.00e-2 & 1.4120 & & 4.81e-3 & 1.1055 & & 8.45e-5 & 3.8546 & & 7.96e-6 & 3.9945 \\  
160 & 8.84e-3 & 1.7618 & & 2.07e-3 & 1.2182 & & 3.38e-6 & 4.6451 & & 4.98e-7 & 3.9989 \\ 
\hline
\end{tabular}
\label{tab:advection_1d_Linf}
\end{table}

\subsection{1D Euler equations}
Next, we solve the one-dimensional Euler equations of gas dynamics
\begin{equation} \label{eq:euler_1d}
 \bfu_t + \bff(\bfu)_x = 0, 
\end{equation}
where the column vector $\bfu$ of the conserved variables and the flux vector $\bff$ are given by 
$$
   \bfu = \left[ \rho,~\rho u,~E \right]^T, \quad \bff(\bfu) = \left[ \rho u,~\rho u^2 + P,~u(E+P) \right]^T.
$$
with $\rho,\: u,\: P$ and $E$ denoting density, velocity, pressure and the specific kinetic energy, respectively.

\begin{figure}[htbp]
\centering
\includegraphics[width=0.45\textwidth]{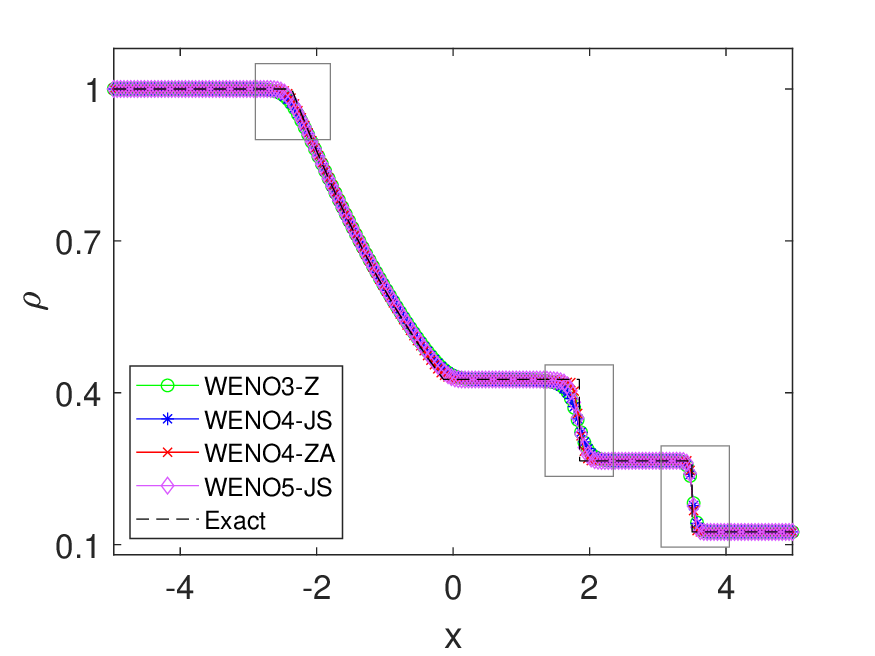}
\includegraphics[width=0.45\textwidth]{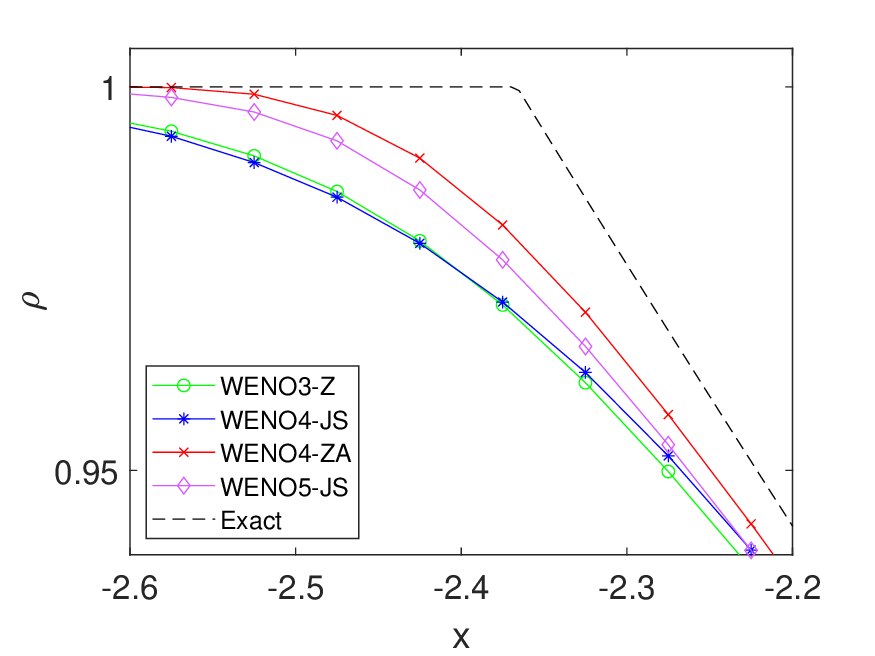}
\includegraphics[width=0.45\textwidth]{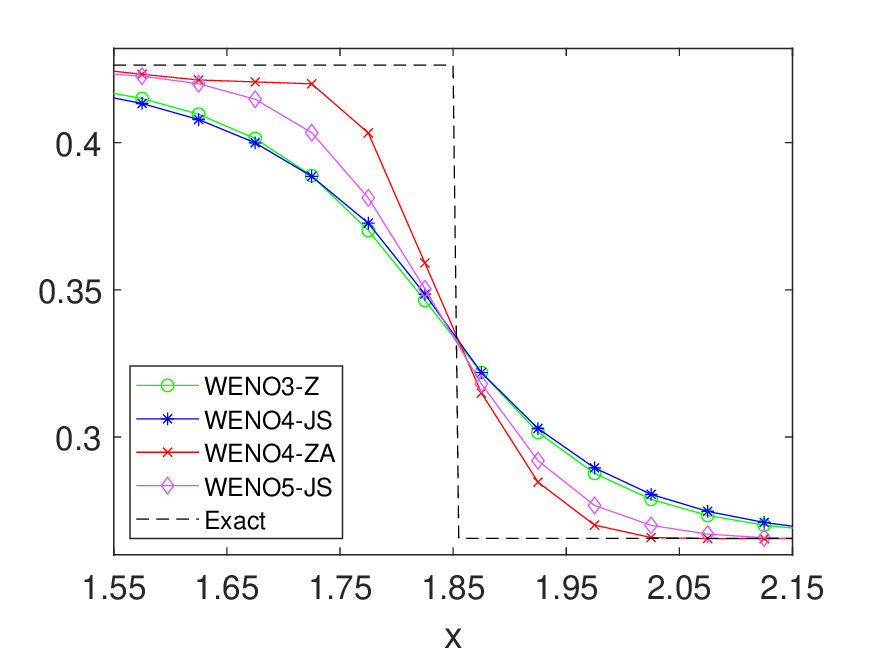}
\includegraphics[width=0.45\textwidth]{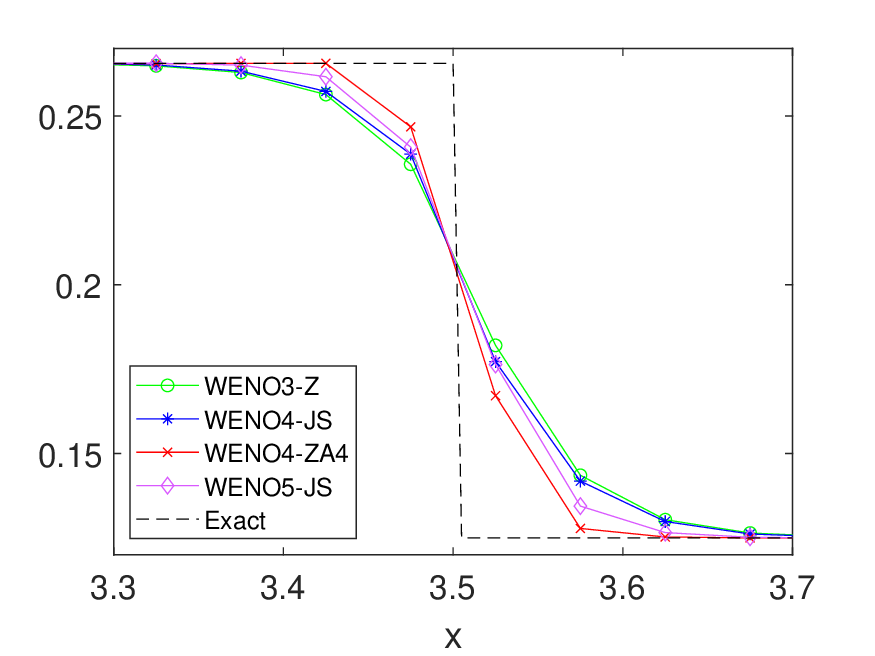}
\caption{Density profiles for the Sod's problem \eqref{eq:euler_1d} and \eqref{eq:sod} at $T=2$ (top left), close-up view of the solutions in the boxes from left to right (top right, bottom left, bottom right) approximated by WENO3-Z (green), WENO4-JS (blue), WENO4-ZA (red) and WENO5-JS (purple) with $N = 200$. 
The dashed black lines are the exact solution.}
\label{fig:sod}
\end{figure}

\begin{figure}[htbp]
\centering
\includegraphics[width=0.45\textwidth]{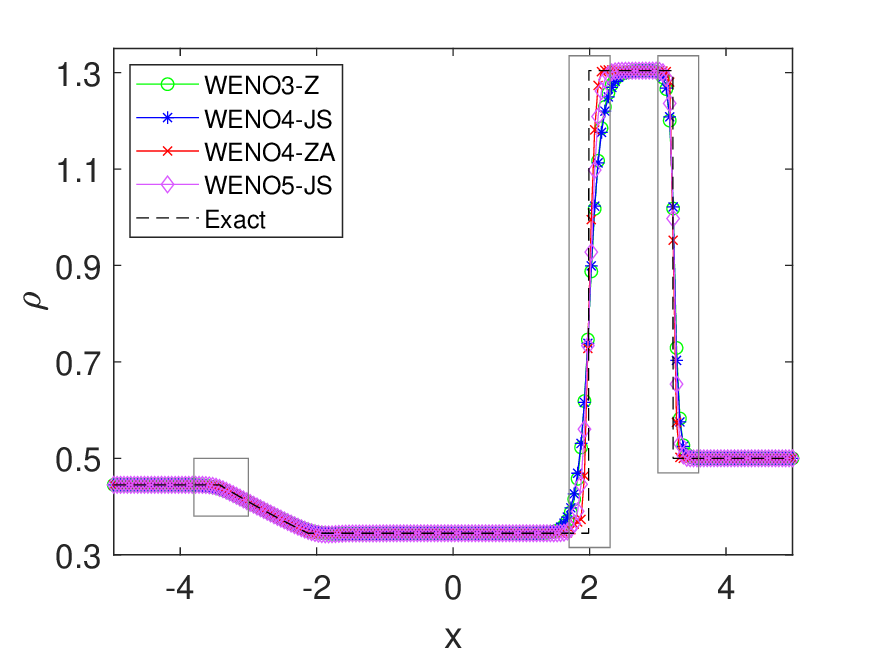}
\includegraphics[width=0.45\textwidth]{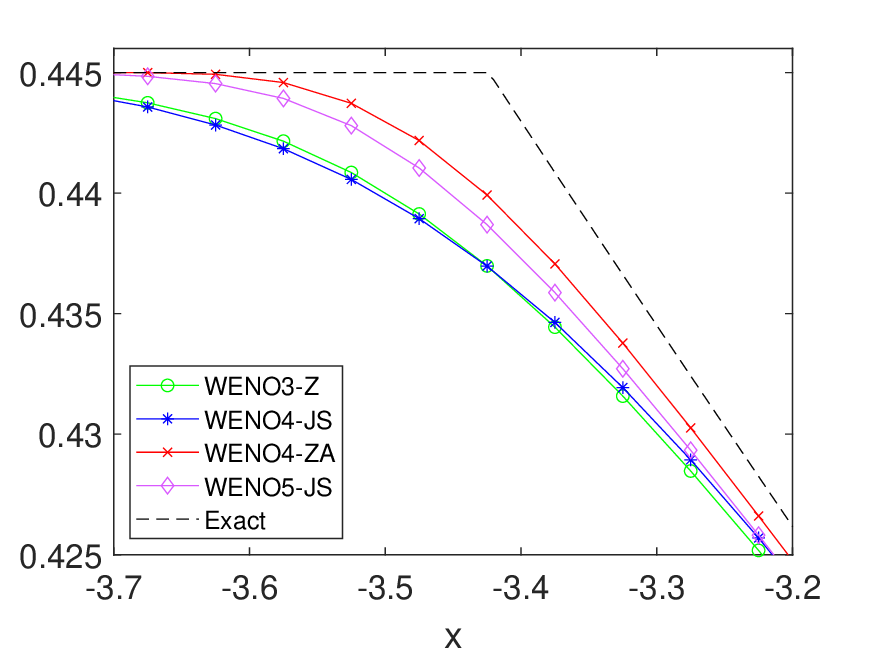}
\includegraphics[width=0.45\textwidth]{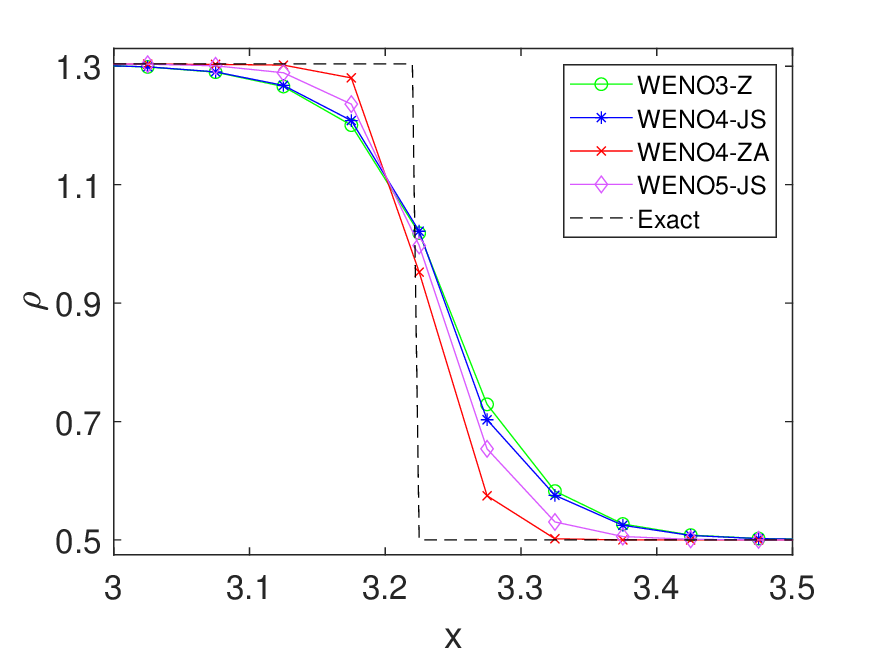}
\includegraphics[width=0.45\textwidth]{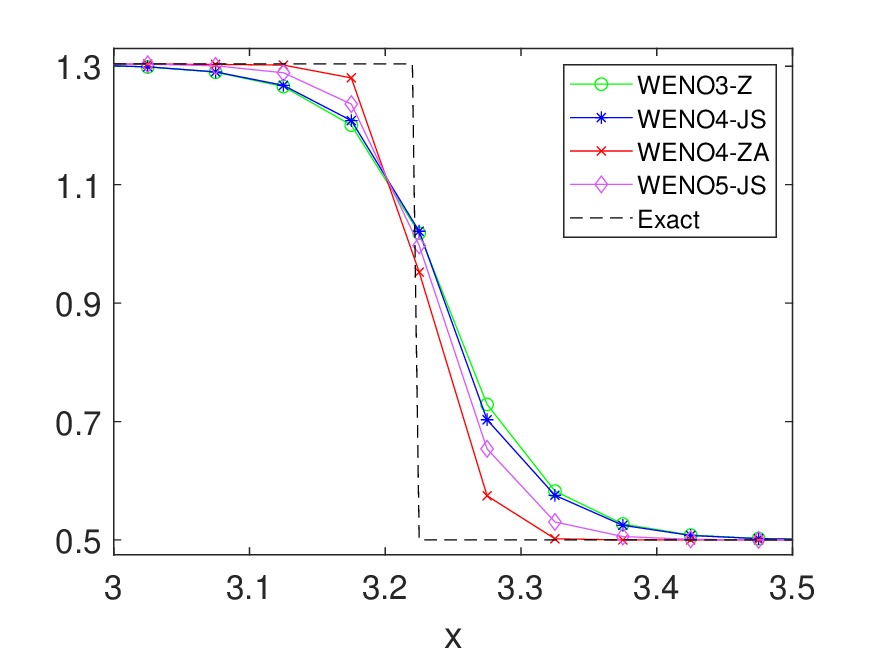}
\caption{Density profiles for the Lax's problem \eqref{eq:euler_1d} and \eqref{eq:lax} at $T=1.3$ (top left), close-up view of the solutions in the boxes from left to right (top right, bottom left, bottom right) solved by WENO3-Z (green), WENO4-JS (blue), WENO4-ZA (red) and WENO5-JS (purple) with $N = 200$.
The dashed black lines are the exact solution.}
\label{fig:lax}
\end{figure}

\begin{example} \label{ex:shock_tube}
We employ the shock tube problem with two sets of the Riemann-type initial data to test the performance of the proposed WENO4-JS  and WENO4-ZA ($p=100$) schemes.
The computational domain $[-5,\: 5]$ is divided into $N = 200$ uniform cells.
For all figures in this example, the top left figure plots the numerical solutions on the entire domain while the other three figures show the solution profiles in the boxes specified in the top left figure.

Consider the Sod's problem with the initial condition
\begin{equation} \label{eq:sod}
 (\rho, u, P ) = \left\{ 
                  \begin{array}{ll} 
                   (1,~0,~1),       & x \leqslant 0, \\ 
                   (0.125,~0,~0.1), & x > 0.
                  \end{array} 
                 \right.
\end{equation}
We integrate the equation until the final time $T = 2$.
In Figure \ref{fig:sod} we show the approximate density $\rho$ at the final time.
The density profile produced by WENO4-JS is comparable to WENO3-Z.
The density approximation by WENO4-ZA shows the sharper solution around the regions of rarefaction, contact discontinuity and shock wave than WENO5-JS. 
The observation is consistent with the $L_1,\: L_2$ and $L_\infty$ errors in Table \ref{tab:sod}, where WENO4-ZA performs better than the other three WENO schemes in terms of accuracy.

For the Lax's problem, the initial condition is
\begin{equation} \label{eq:lax}
 (\rho, u, P ) = \left\{ 
                  \begin{array}{ll} 
                   (0.445,~0.698,~3.528), & x \leqslant 0, \\ 
                   (0.5,~0,~0.571),       & x > 0. 
                  \end{array} 
                 \right. 
\end{equation}
Figure \ref{fig:lax} shows the approximate density from each WENO scheme at the final time $T = 1.3$. 
As in the Sod's problem, the density approximated by WENO4-JS looks similar to the one by WENO3-Z, while the density profile with WENO4-ZA exhibits less dissipation near the rarefied, contact discontinuity and shock regions than WENO5-JS.
Table \ref{tab:lax} shows that WENO3-Z, WENO4-JS, WENO5-JS and WENO4-ZA yield increasingly more accurate approximations.
\end{example}

\begin{table}[htbp]
\renewcommand{\arraystretch}{1.1}
\scriptsize
\centering
\caption{$L_1,\: L_2,\: L_\infty$ errors for the Sod's problem \eqref{eq:euler_1d} and \eqref{eq:sod}.}      
\begin{tabular}{ccccc}
\hline 
Error & WENO3-Z & WENO4-JS & WENO4-ZA & WENO5-JS \\
\hline
$L_1$      & 4.958e-3 & 4.834e-3 & 2.323e-3 & 3.476e-3 \\
$L_2$      & 1.176e-2 & 1.152e-2 & 7.590e-3 & 9.673e-3 \\
$L_\infty$ & 7.989e-2 & 7.778e-2 & 6.717e-2 & 7.595e-2 \\  
\hline
\end{tabular}
\label{tab:sod}
\end{table}

\begin{table}[htbp]
\renewcommand{\arraystretch}{1.1}
\scriptsize
\centering
\caption{$L_1,\: L_2,\: L_\infty$ errors for the Lax's problem \eqref{eq:euler_1d} and \eqref{eq:lax}.}      
\begin{tabular}{ccccc}
\hline
Error & WENO3-Z & WENO4-JS & WENO4-ZA & WENO5-JS \\
\hline
$L_1$      & 1.753e-2 & 1.770e-2 & 8.334e-3 & 1.203e-2 \\
$L_2$      & 6.870e-2 & 6.815e-2 & 4.926e-2 & 5.846e-2 \\
$L_\infty$ & 5.180e-1 & 5.214e-1 & 4.523e-1 & 4.973e-1 \\  
\hline
\end{tabular}
\label{tab:lax}
\end{table}

\begin{figure}[htbp]
\centering
\includegraphics[width=0.45\textwidth]{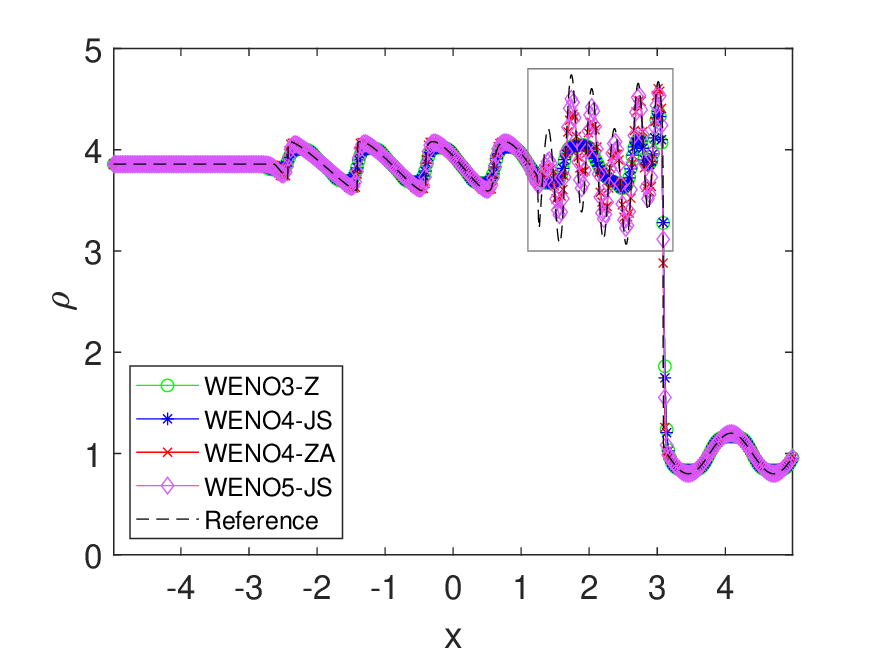}
\includegraphics[width=0.45\textwidth]{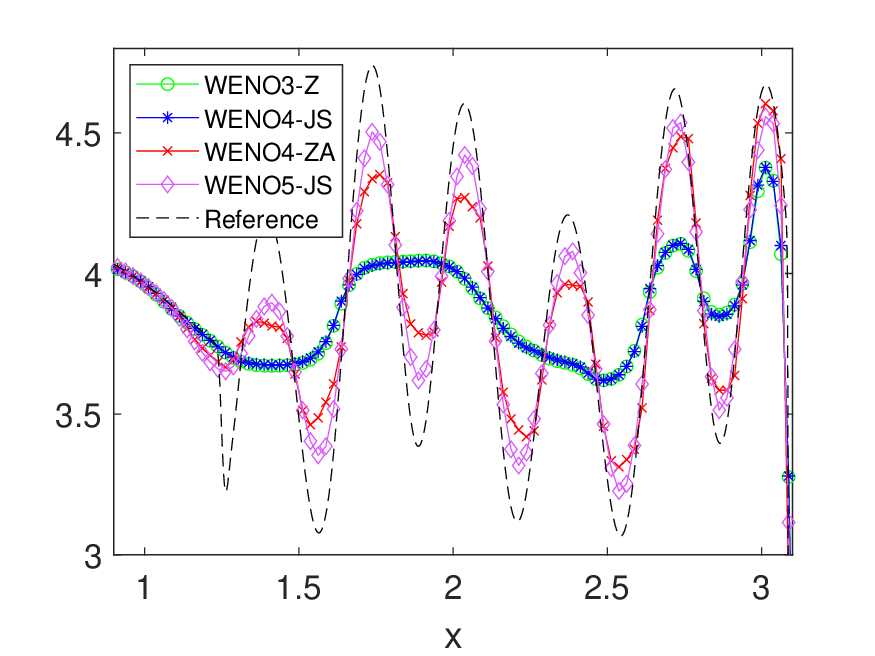}
\caption{Solution profiles for Example \ref{ex:shock_entropy_wave} with $k=5$ at $T=2$ (left), close-up view of the solutions in the box (right) computed by WENO3-Z (green), WENO4-JS (blue), WENO4-ZA (red) and WENO5-JS (purple) with $N=400$.
The dashed black lines are generated by WENO5-M with $N=4000$.}
\label{fig:shock_entropy_wave_k5}
\end{figure}

\begin{figure}[htbp]
\centering
\includegraphics[width=0.45\textwidth]{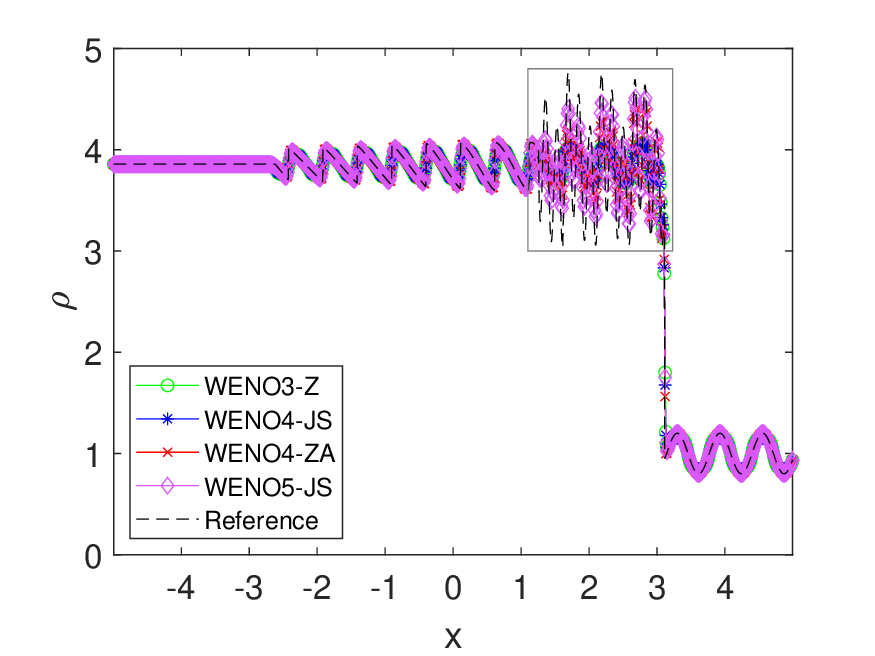}
\includegraphics[width=0.45\textwidth]{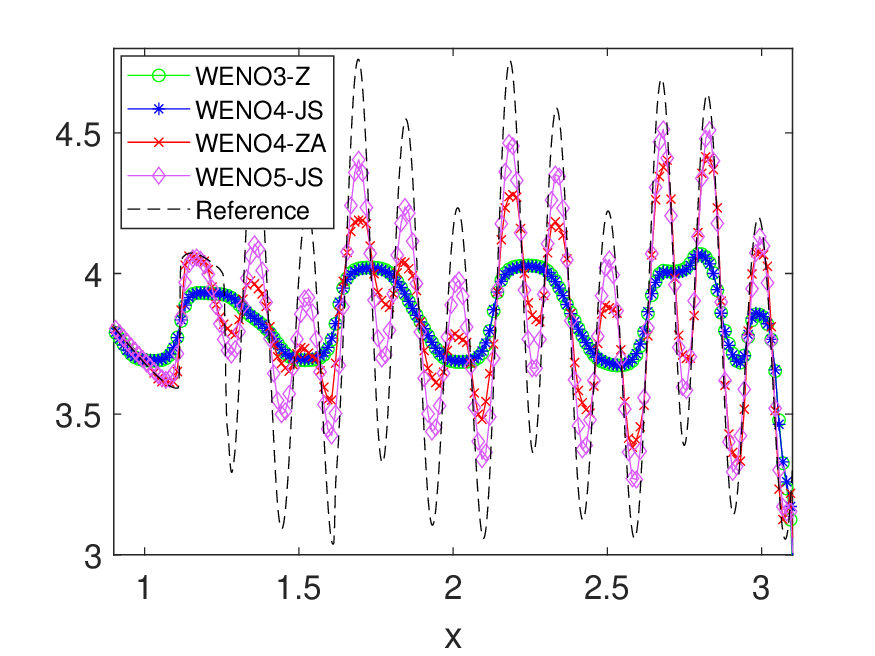}
\caption{Solution profiles for Example \ref{ex:shock_entropy_wave} with $k=10$ at $T=2$ (left), close-up view of the solutions in the box (right) computed by WENO3-Z (green), WENO4-JS (blue), WENO4-ZA (red) and WENO5-JS (purple) with $N=800$.
The dashed black lines are generated by WENO5-M with $N=8000$.}
\label{fig:shock_entropy_wave_k10}
\end{figure}

\begin{example} \label{ex:shock_entropy_wave}
The shock density wave interaction problem \cite{ShuOsherII}, which handles a right moving Mach 3 shock interacting with sine waves in density, has the initial condition
$$
   (\rho, u, P ) = \left\{ 
                    \begin{array}{ll} 
                     (3.857143,~2.629369,~10.333333), & x < -4, \\ 
                     (1 + 0.2 \sin(kx),~0,~1),        & x \geqslant -4,
                    \end{array} 
                   \right. 
$$
where $k$ is the wave number of the entropy wave.
We pick $p=100$ for WENO-ZA.
For $k=5$, we divide the computational domain $[-5, 5]$ into $N=400$ uniform cells.
Figure \ref{fig:shock_entropy_wave_k5} plots the density profile approximated by those WENO schemes at $T = 2$.
For $k=10$, the domain $[-5, 5]$ is divided into $N = 800$ uniform cells.
We present the numerical solutions of the density at $T=2$ in Figure \ref{fig:shock_entropy_wave_k10}.
We observe that the simulations from WENO3-Z and WENO4-JS are comparable, and WENO4-JS performs better in capturing the fine structure in the density profile than WENO3-Z and WENO4-JS, even if WENO4-ZA produces slightly less accurate result than WENO5-JS.
\end{example}

\begin{figure}[htbp]
\centering
\includegraphics[width=0.325\textwidth]{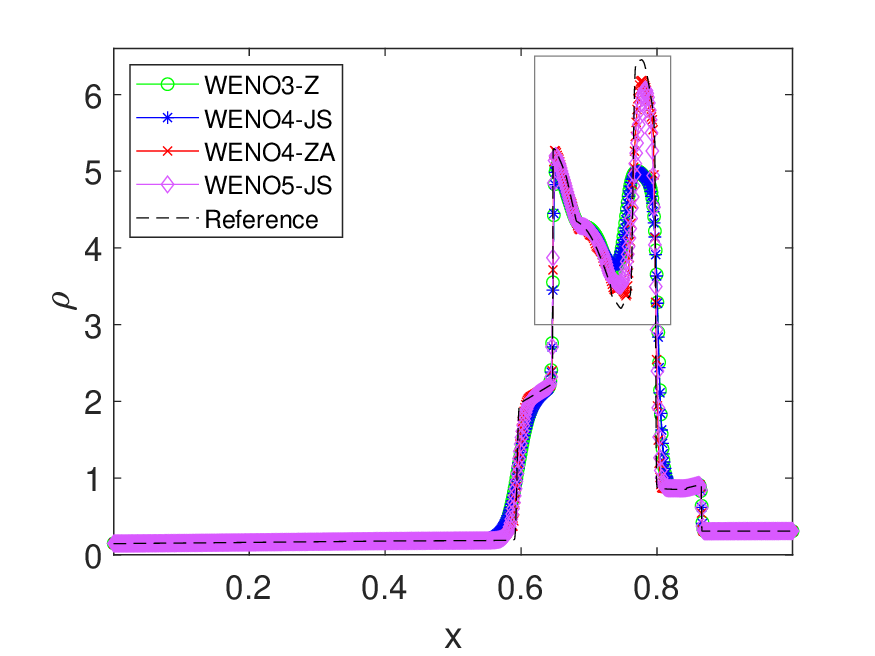}
\includegraphics[width=0.325\textwidth]{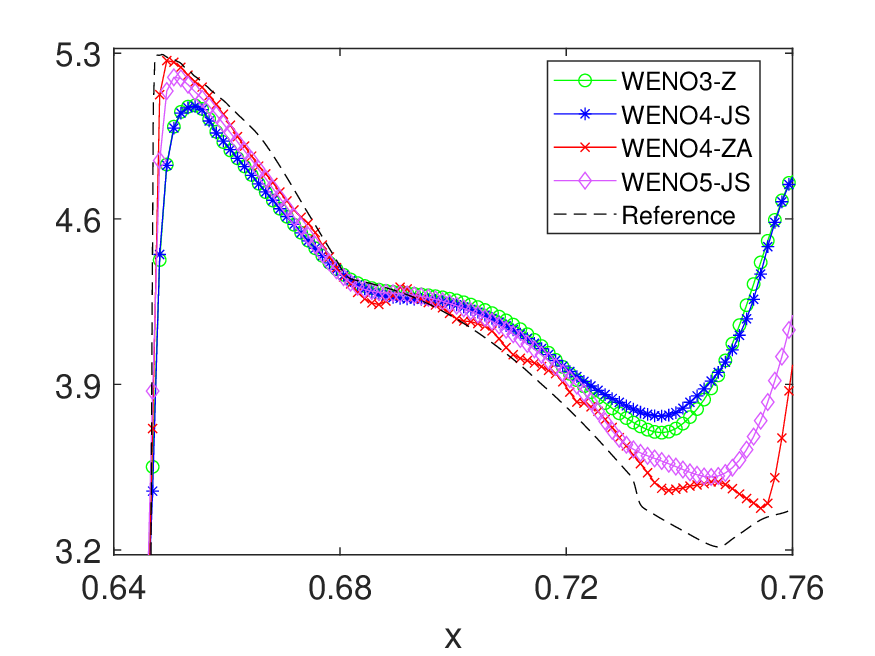}
\includegraphics[width=0.325\textwidth]{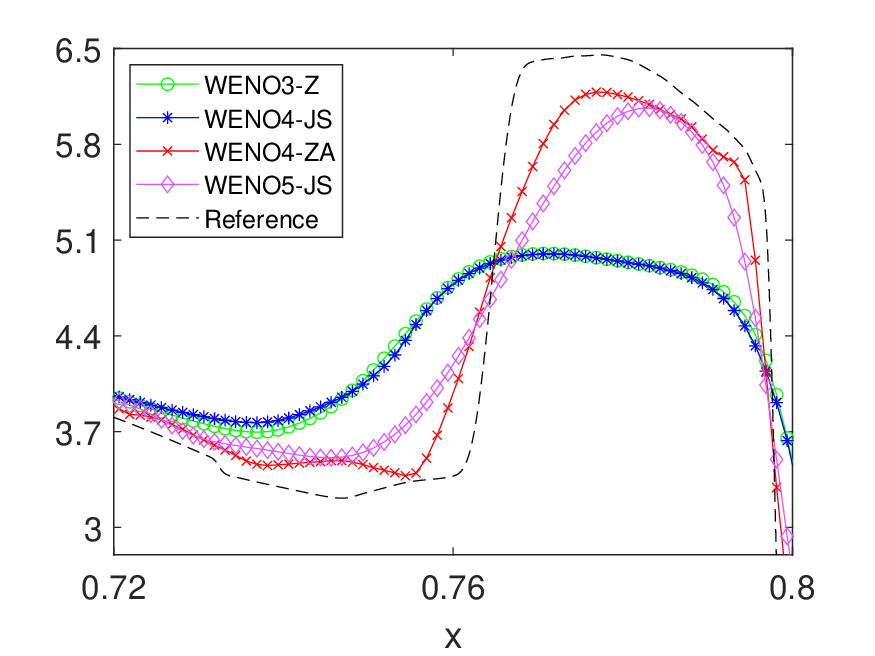}
\caption{Density profiles for Example \ref{ex:interacting_blastwave} at $T=0.038$ (left), close-up views of the solutions in the box (middle, right) computed by WENO3-Z (green), WENO4-JS (blue), WENO4-ZA (red) and WENO5-JS (purple) with $N = 800$.
The dashed black lines are generated by WENO5-M with $N = 4000$.}
\label{fig:interacting_blastwave}
\end{figure}

\begin{example} \label{ex:interacting_blastwave}
Consider the blastwaves interaction problem \cite{Woodward} with the initial condition
$$
   (\rho, u, P ) = \left\{ 
                    \begin{array}{ll} 
                     (1,~0,~1000), & 0 \leqslant x < 0.1, \\ 
                     (1,~0,~0.01), & 0.1 \leqslant x < 0.9, \\
                     (1,~0,~100),  & 0.9 \leqslant x \leqslant 1.
                    \end{array} 
                   \right. 
$$
We divide the computational domain $[0, 1]$ into $N = 800$ uniform cells.
The numerical solutions of the density $\rho$ at $T = 0.038$ are shown in Figure \ref{fig:interacting_blastwave}.
The density profile produced by WENO4-JS is comparable to WENO3-Z.
In the close-up regions, WENO4-ZA ($p=100$) has better resolution than WENO5-JS.
\end{example}

\subsection{2D Euler equations} 
The two-dimensional Euler equations of gas dynamics have the form
\begin{equation} \label{eq:euler_2d}
 \bfu_t + \bff(\bfu)_x + \bfg(\bfu)_y = 0, 
\end{equation}
where the conserved vector $\bfu$ and the flux functions $\bff, \bfg$ in the respective $x, y$ directions are 
\begin{align*}
      \bfu  &= \left[ \rho, ~\rho u, ~\rho v, ~E \right]^T, \\
 \bff(\bfu) &= \left[ \rho u, ~\rho u^2 + P, ~\rho u v, ~u(E+P) \right]^T, \\ 
 \bfg(\bfu) &= \left[ \rho v, ~\rho u v, ~\rho v^2 + P, ~v(E+P) \right]^T.
\end{align*}
As in one-dimensional case, $\rho$ is the density, $P$ is the pressure and $E$ is the specific kinetic energy. 
The primitive variables $u$ and $v$ denote $x$- and $y$-component velocities, respectively.

\begin{figure}[h!]
\centering
\includegraphics[width=0.45\textwidth]{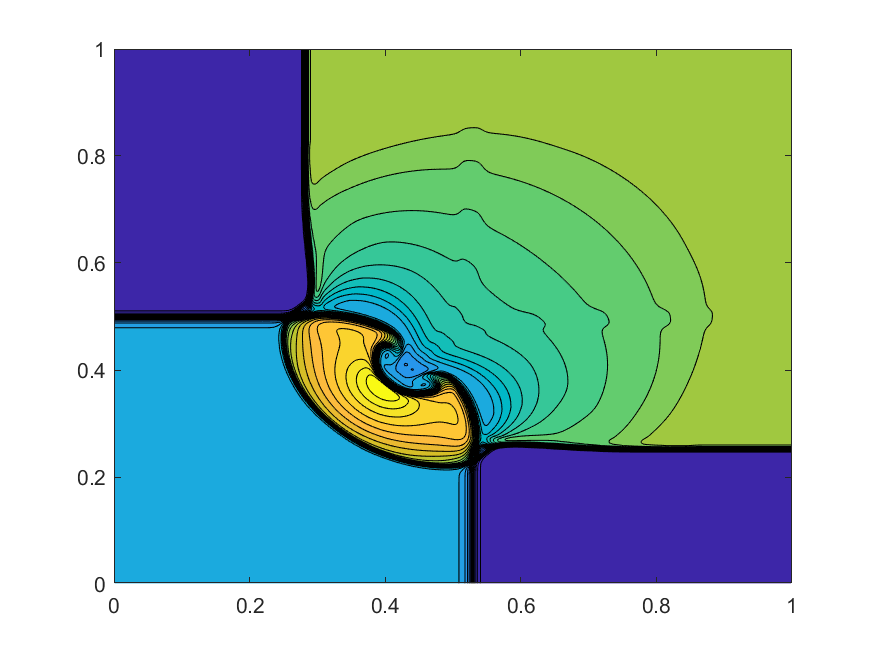}
\includegraphics[width=0.45\textwidth]{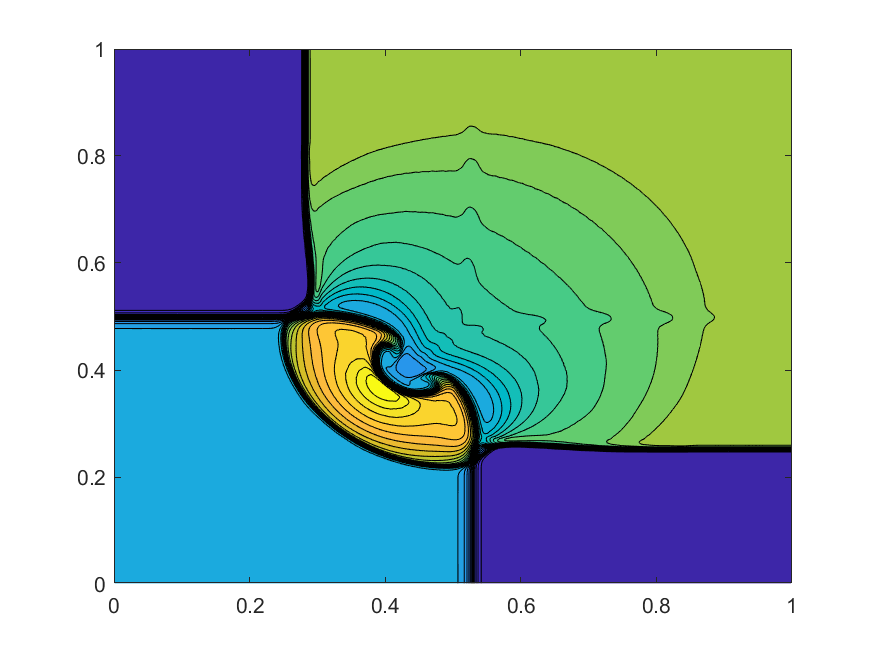}
\includegraphics[width=0.45\textwidth]{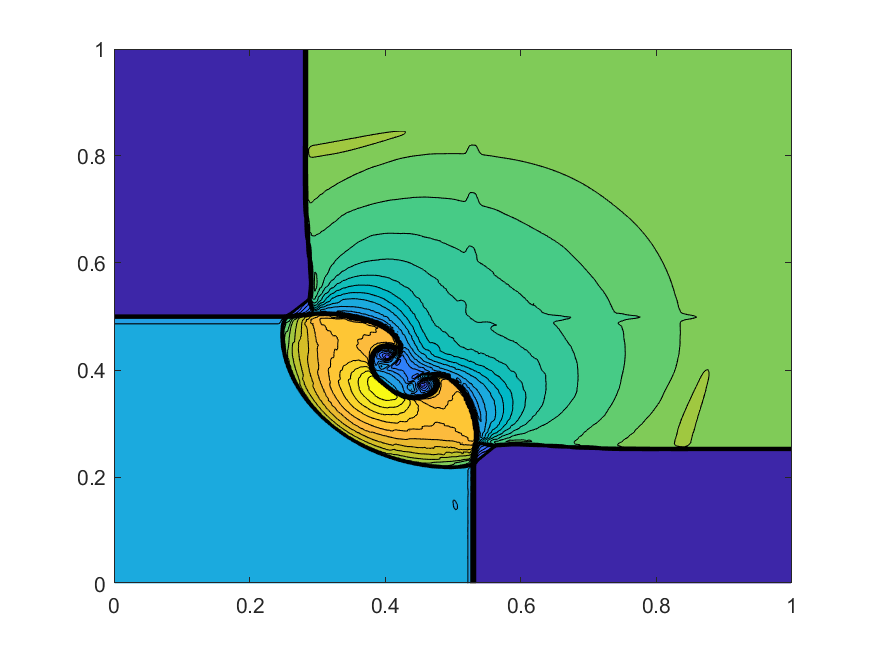}
\includegraphics[width=0.45\textwidth]{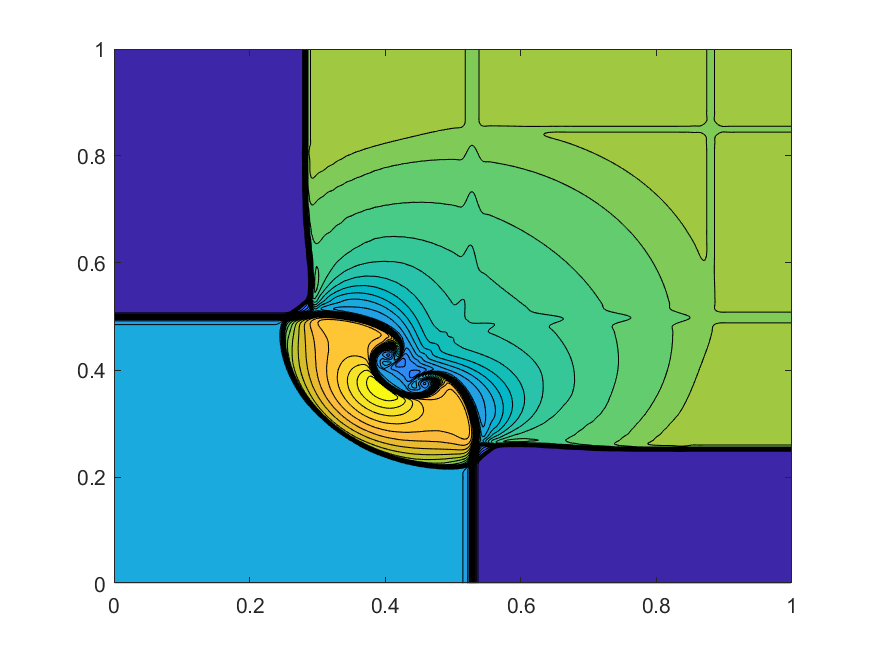}
\caption{Density in the filled contour plot for Example \ref{ex:euler_2d_riemann} at $T=0.3$ by WENO3-Z (top left), WENO4-JS (top right), WENO4-ZA (bottom left) and WENO5-JS (bottom right) with $N_x = N_y = 400$.
Each contour plot displays contours at $30$ levels of the density.}
\label{fig:eer2d}
\end{figure}

\begin{example} \label{ex:euler_2d_riemann}
We consider the Riemann problem in \cite{Kurganov} for the two-dimensional Euler equations \eqref{eq:euler_2d}, along with the initial condition,
$$
   (\rho, u, v, P ) = \left\{ 
                       \begin{array}{ll} 
                        (1,~0.1,~0,~1),             & x > 0.5,~y > 0.5, \\
                        (0.5313,~0.8276,~0,~0.4),   & x < 0.5,~y > 0.5, \\
                        (0.8,~0.1,~0,~0.4),         & x < 0.5,~y < 0.5, \\
                        (0.5313,~0.1,~0.7276,~0.4), & x > 0.5,~y < 0.5,
                       \end{array}
                      \right.
$$
We divide the square computational domain $[0,1] \times [0,1]$ into $N_x \times N_y = 400 \times 400$ uniform cells.
The free stream boundary conditions are imposed in both $x$- and $y$-directions.
The simulations of the density computed by WENO3-Z, WENO4-JS, WENO4-ZA ($p=100$) and WENO5-JS at the final time $T=0.3$ are presented in Figure \ref{fig:eer2d}.
It is observed that all WENO schemes capture the ripples well, which is in agreement with \cite{Kurganov}.
\end{example}

\begin{figure}[h!]
\centering
\includegraphics[width=0.48\textwidth]{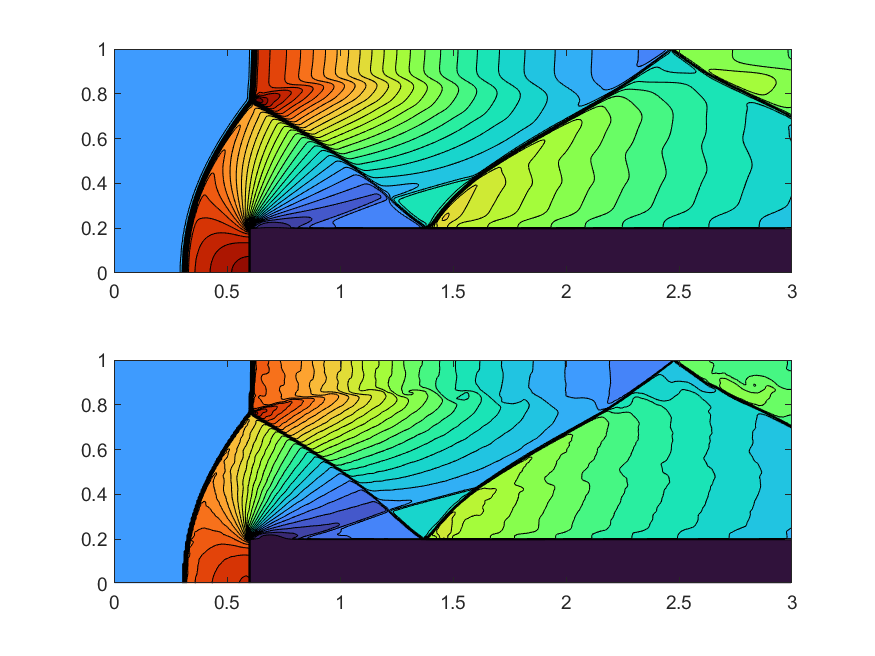}
\includegraphics[width=0.48\textwidth]{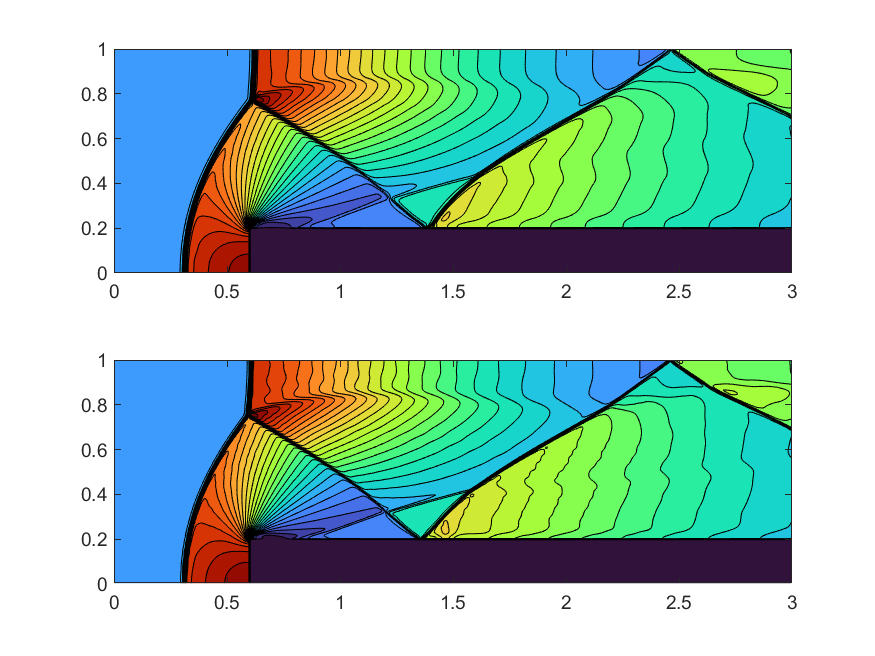}
\caption{Density in the filled contour plot for Example \ref{ex:forward_facing_step} at $T=4$ by WENO3-Z (top left), WENO4-JS (top right), WENO4-ZA (bottom left) and WENO5-JS (bottom right) with $N_x = 480$ and $N_y = 160$.
Each contour plot displays contours at $30$ levels of the density.}
\label{fig:ffs}
\end{figure}

\begin{example} \label{ex:forward_facing_step}
The forward facing step problem \cite{Woodward} involves a Mach 3 flow at the inlet to a rectangular wind tunnel with a step near the inlet region that creates a transient shock wave pattern.
Initially, the tunnel is full of the gas with $(\rho, u, v, P ) = (1.4, 3, 0, 1)$.
The computational domain $[0,3] \times [0,1]$ is divided into $N_x \times N_y = 480 \times 160$ uniform cells.
Inflow boundary condition is imposed at the inlet, and all of the gradients vanish at the outlet.
The reflective boundary conditions are applied along the walls of the tunnel. 
We compute the numerical solutions until the final time $T=4$ and present the contour lines of the approximate density in Figure \ref{fig:ffs}.
We notice that the vortical structures emulated along the slip line are less diffused with WENO4-ZA ($p=20$) and WENO5-JS than WENO3-Z and WENO4-JS.
\end{example}

\begin{figure}[h!]
\centering
\includegraphics[width=0.45\textwidth]{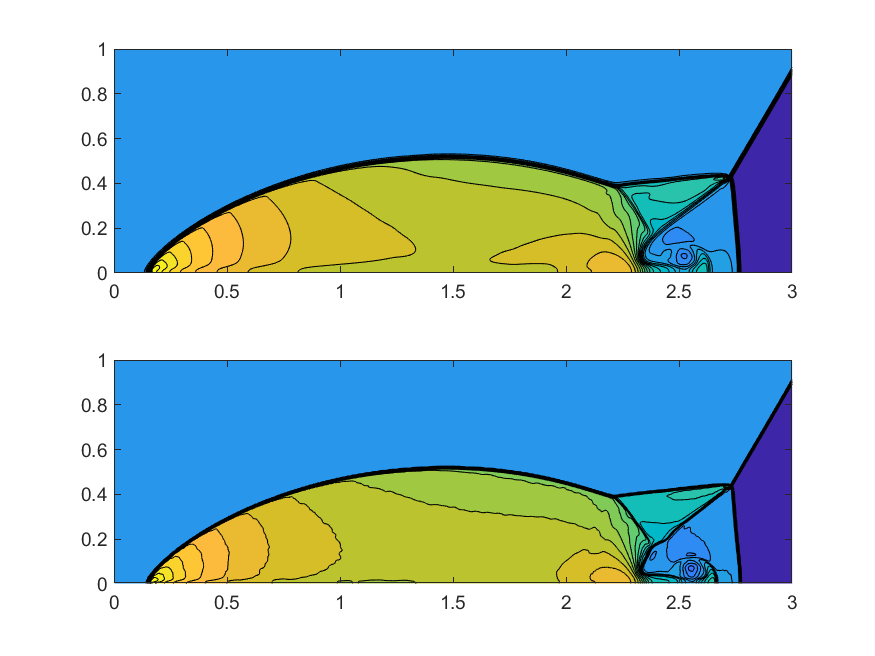}
\includegraphics[width=0.45\textwidth]{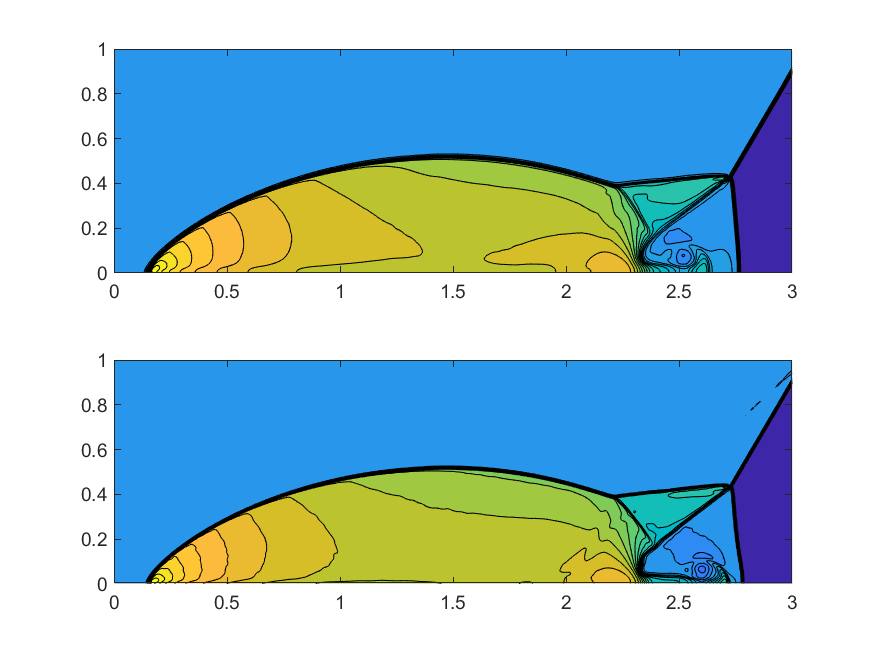}
\caption{Density in the filled contour plot for Example \ref{ex:double_mach_reflection} at $T=0.2$ by WENO3-Z (top left), WENO4-JS (top right), WENO4-ZA (bottom left) and WENO5-JS (bottom right) with $N_x = 800$ and $N_y = 200$.
Each contour plot displays contours at $30$ levels of the density.}
\label{fig:dmr}
\end{figure}

\begin{example} \label{ex:double_mach_reflection}
We end by solving the double Mach reflection problem \cite{Woodward}, which describes the reflection of a planar Mach shock in the air hitting a wedge.
The initial condition is given by 
$$
   (\rho, u, v, P ) = \left\{ 
                       \begin{array}{ll} 
                        (8,~8.25 \cos \theta,~-8.25 \sin \theta,~116.5), & x < \frac{1}{6} + \frac{y}{\sqrt{3}}, \\ 
                        (1.4,~0,~0,~1), & x \geqslant \frac{1}{6} + \frac{y}{\sqrt{3}}, 
                       \end{array} 
                      \right. 
$$
with $\theta = \frac{\pi}{6}$.
The computational domain is $[0,4] \times [0,1]$ with $N_x \times N_y = 800 \times 200$ uniform cells. 
For the bottom boundary, the exact post-shock condition is imposed for $\left[ 0, 1/6 \right)$ on the $x$-axis whereas the reflective boundary condition is applied to the rest.
At the top boundary, the exact motion of the right-moving Mach $10$ oblique shock is used.
The supersonic inflow and free stream outflow boundary conditions are imposed for the left and right boundary, respectively.
The simulation is conducted until the final time $T = 0.2$, and the numerical densities are shown in Figure \ref{fig:dmr}.
We can see that the small-scale structures along the slip line and the vortical rollup at the tip of the jet are better resolved in WENO4-ZA ($p=10$) and WENO5-JS.
\end{example}

\section{Conclusion} \label{sec:conclusion}
In this study, we have developed novel JS-type and Z-type nonlinear weights, based on simple local smoothness indicators, for the fourth-order central-upwind WENO scheme.
The sufficient condition for fourth-order accuracy is derived, where the JS-type weights do not satisfy while the Z-type weights do if we properly set the value of $p$ in \eqref{eq:weights4_ZA}.
We also analyze the weights behavior around discontinuity and find that the Z-type weights are closer to the linear weights than the JS-type ones, where the former nonlinear weights are expected to yield less dissipation than the latter ones.
The accuracy test show that the WENO4-ZA scheme achieves fourth-order accuracy.
The one- and two-dimensional benchmarks of Euler equations demonstrate that the proposed central-upwind WENO schemes work well for contact waves and fine-scale structures.

\section*{Acknowledgments}
We would like to thank Baoshan Wang for helpful discussions.
This work is supported by National Research Foundation of Korea (NRF) under the grant number 2021R1A2C3009648, by POSTECH Basic Science Research Institute under the NRF grant number 2021R1A6A1A10042944, and partially by NRF grant funded by the Korea government (MSIT) (No. RS-2023-00219980). 
This work was supported by POSTECH Basic Science Research Institute under the NRF grant number NRF2021R1A6A1A1004294412.


\begin{thebibliography}{10}

\bibitem{Borges}
Rafael Borges, Monique Carmona, Bruno Costa, and Wai~Sun Don, \emph{An improved
  weighted essentially non-oscillatory scheme for hyperbolic conservation
  laws}, J. Comput. Phys. \textbf{227} (2008), no.~6, 3191--3211.

\bibitem{ChenGuJung}
Xinjuan Chen, Jiaxi Gu, and Jae-Hun Jung, \emph{A spatial-temporal weight
  analysis and novel nonlinear weights of weighted essentially non-oscillatory
  schemes for hyperbolic conservation laws}, J. Sci. Comput. \textbf{102}
  (2025), 34.

\bibitem{Don}
Wai-Sun Don and Rafael Borges, \emph{Accuracy of the weighted essentially
  non-oscillatory conservative finite difference schemes}, J. Comput. Phys.
  \textbf{250} (2013), 347--372.

\bibitem{Gu}
Jiaxi Gu, Xinjuan Chen, and Jae-Hun Jung, \emph{Fifth-order weighted
  essentially non-oscillatory schemes with new {Z}-type nonlinear weights for
  hyperbolic conservation laws}, Comput. Math. Appl. \textbf{134} (2023),
  140--166.

\bibitem{Henrick}
Andrew~K. Henrick, Tariq~D. Aslam, and Joseph~M. Powers, \emph{Mapped weighted
  essentially non-oscillatory schemes: Achieving optimal order near critical
  points}, J. Comput. Phys. \textbf{207} (2005), no.~2, 542--567.

\bibitem{Hu}
Fuxing Hu, \emph{The 6th-order weighted {ENO} schemes for hyperbolic
  conservation laws}, Comput. \& Fluids \textbf{174} (2018), 34--45.

\bibitem{HuAdams}
X.Y. Hu and N.A. Adams, \emph{Scale separation for implicit large eddy
  simulation}, J. Comput. Phys. \textbf{230} (2011), no.~19, 7240--7249.

\bibitem{HuWangAdams}
X.Y. Hu, Q.~Wang, and N.A. Adams, \emph{An adaptive central-upwind weighted
  essentially non-oscillatory scheme}, J. Comput. Phys. \textbf{229} (2010),
  no.~23, 8952--8965.

\bibitem{Huang}
Cong Huang and Li~Li Chen, \emph{A new adaptively central-upwind sixth-order
  weno scheme}, J. Comput. Phys. \textbf{357} (2018), 1--15.

\bibitem{Jiang}
Guang-Shan Jiang and Chi-Wang Shu, \emph{Efficient implementation of weighted
  {ENO} schemes}, J. Comput. Phys. \textbf{126} (1996), no.~1, 202--228.

\bibitem{Kurganov}
Alexander Kurganov and Eitan Tadmor, \emph{Solution of two-dimensional
  {Riemann} problems for gas dynamics without {Riemann} problem solvers},
  Numer. Methods Partial Differential Equations \textbf{18} (2002), no.~5,
  584--608.

\bibitem{LiuShenPengZhang}
Shengping Liu, Yiqing Shen, Jun Peng, and Jun Zhang, \emph{Two-step weighting
  method for constructing fourth-order hybrid central {WENO} scheme}, Comput.
  \& Fluids \textbf{207} (2020), 104590.

\bibitem{ShuSpringer}
Chi-Wang Shu, \emph{Essentially non-oscillatory and weighted essentially
  non-oscillatory schemes for hyperbolic conservation laws}, Advanced Numerical
  Approximation of Nonlinear Hyperbolic Equations. Lecture Notes in Mathematics
  (A.~Quarteroni, ed.), Springer, Berlin, 1998, pp.~325--432.

\bibitem{ShuOsherI}
Chi-Wang Shu and Stanley Osher, \emph{Efficient implementation of essentially
  non-oscillatory shock-capturing schemes}, J. Comput. Phys. \textbf{77}
  (1988), no.~2, 439--471.

\bibitem{ShuOsherII}
\bysame, \emph{Efficient implementation of essentially non-oscillatory
  shock-capturing schemes, {II}}, J. Comput. Phys. \textbf{83} (1989), no.~1,
  32--78.

\bibitem{WangDonLiWang}
Cai-Feng Wang, Wai~Sun Don, Jia-Le Li, and Bao-Shan Wang, \emph{Improved
  sixth-order {WENO} finite difference schemes for hyperbolic conservation
  laws}, Adv. Appl. Math. Mech. (2025), ,accepted for publication.

\bibitem{Wang}
Yahui Wang, Yulong Du, Kunlei Zhao, and Li~Yuan, \emph{A new 6th-order {WENO}
  scheme with modified stencils}, Comput. \& Fluids \textbf{208} (2020),
  104625.

\bibitem{Woodward}
Paul Woodward and Phillip Colella, \emph{The numerical simulation of
  two-dimensional fluid flow with strong shocks}, J. Comput. Phys. \textbf{54}
  (1984), no.~1, 115--173.

\bibitem{Zhao}
Guo‑Yan Zhao, Ming‑Bo Sun, Yong Mei, Liang Li, Hong‑Bo Wang, Guang‑Xin
  Li, Yuan Liu, Yong‑Chao Sun, and Chang‑Hai Liang, \emph{An efficient
  adaptive central-upwind {WENO-CU6} numerical scheme with a new sensor}, J.
  Sci. Comput. \textbf{81} (2019), 649–670.

\end{thebibliography}

\providecommand{\bysame}{\leavevmode\hbox to3em{\hrulefill}\thinspace}
\providecommand{\MR}{\relax\ifhmode\unskip\space\fi MR }
\providecommand{\MRhref}[2]{%
  \href{http://www.ams.org/mathscinet-getitem?mr=#1}{#2}
}
\providecommand{\href}[2]{#2}

\end{document}